\documentclass[11pt,a4paper]{amsart}

\usepackage{lineno}

\usepackage[latin1]{inputenc}
\usepackage[english]{babel}
\usepackage{amssymb,amsfonts,amsmath,mathtools}

\usepackage{graphicx}
\usepackage{braket}
\usepackage[pdftex,plainpages=false,colorlinks,hyperindex,bookmarksopen,linkcolor=red,citecolor=blue,urlcolor=blue]{hyperref}

\usepackage{cite}

\usepackage{bm}

\usepackage{mathrsfs}

\usepackage{epstopdf}

\usepackage{braket}

\usepackage[hmargin=3cm,vmargin={3.5cm,4cm}]{geometry}

\theoremstyle{definition}                                 

\theoremstyle{definition}                           

\theoremstyle{remark}                             

\usepackage{color}

\usepackage{mathtools,slashed}

\newcommand{\be}{\begin{eqnarray}}
\newcommand{\ee}{\end{eqnarray}}

\newcommand{\R}{\mathbb{R}}  
\newcommand{\C}{\mathbb{C}} 
\newcommand{\N}{\mathbb{N}} 
\DeclareMathOperator{\real}{Re}
\def\eu{\ensuremath{\mathrm{e}}}
\def\iu{\ensuremath{\mathrm{i}}}
\def\du{\ensuremath{\mathrm{d}}}
\def\texttiny#1{{\text{\tiny{#1}}}}
\def\DC{{}^{\texttiny{C}}\! D}   
\def\DR{{}^{\texttiny{RL}}\! D}  
\def\JPr{{\bm{\mathcal{J}}}}     
\def\DPrRL{{}^{\texttiny{RL}}\!{\bm{\mathcal{D}}}} 
\def\DPrC{{}^{\texttiny{C}}\!{\bm{\mathcal{D}}}}   

\def\eg{{\it e.g. }} 
\def\ie{{\it i.e. }}

\newcommand{\wt}[1]{\widetilde{#1}}
\newcommand{\wh}[1]{\widehat{#1}}

\newcommand{\ceil}[1]{\lceil #1 \rceil}
\def\d{{\rm d}}

\def\eg{{\it e.g.}\ }
\def\ie{{\it i.e.}\ }

\numberwithin{equation}{section}

\allowdisplaybreaks

\begin{document}

\title{A practical guide to Prabhakar fractional calculus}
		
		\author[A. Giusti, et al.]{Andrea Giusti$^1$}
		\address{${}^1$ 
		Department of Physics $\&$ Astronomy, 
		Bishop's University, 
		2600 College Street, 
		QC J1M 1Z7 Sherbrooke, Canada.}	
 		\email{agiusti@ubishops.ca}
		
		\author[]{Ivano Colombaro$^2$}
		\address{${}^2$ 
		Department of Information and Communication Technologies, 
		Universitat Pompeu Fabra, C/Roc Boronat 138, Barcelona, Spain}	
 		\email{ivano.colombaro@upf.edu}

		\author[]{Roberto Garra$^3$}
		\address{${}^3$ 
		Department of Statistical Sciences, 
		Sapienza University of Rome, Piazzale Aldo Moro 5, 00185 Roma, Italy}	
 		\email{roberto.garra@sbai.uniroma1.it}
		
		\author[]{Roberto Garrappa$^4$}
		\address{${}^4$ 
		Department of Mathematics, 
		University of Bari, Via E. Orabona 4, 70126 Bari, Italy
		and the Italian National Group for Scientific Computing (INdAM-GNCS)}	
 		\email{roberto.garrappa@uniba.it}
		
		\author[]{Federico Polito$^5$}
		\address{${}^5$ 
		Department of Mathematics, 
		University of Torino, 
		Via Carlo Alberto 10, 
		10123 Torino, Italy.}	
 		\email{federico.polito@unito.it}
 		
 		\author[]{Marina Popolizio$^6$}
		\address{${}^6$ 
		Department of Electrical and Information Engineering, 
		Polytechnic University of Bari, 
		Via E. Orabona 4, 
		70126 Bari, Italy
		the Italian National Group for Scientific Computing (INdAM-GNCS)}	
 		\email{marina.popolizio@poliba.it}
		
 		\author[]{Francesco Mainardi$^7$}
		\address{${}^7$ 
		Department of Physics and Astronomy, 
		University of Bologna, Via Irnerio 46, 40126 Bologna, Italy}	
 		\email{mainardi@bo.infn.it}
 				
    \thanks{{\em Fract. Calc. Appl. Anal.}, Vol. {\bf 23}, No 1 (2020), pp. 9--54 ; {\bf doi:} 10.1515/fca-2020-0002}

\begin{abstract}
The Mittag--Leffler function is universally acclaimed as the Queen function of fractional calculus. 
The aim of this work is to survey 
 the key results and applications emerging from the three-parameter generalization of this function, known as the Prabhakar function.
Specifically, after reviewing key historical events that led to the discovery and modern development of this peculiar function, 
we discuss how the latter allows one to introduce an enhanced scheme for fractional calculus. 
Then, we summarize the progress in the application of this new general framework to physics and renewal processes. 
We also provide a collection of results on the numerical evaluation of the Prabhakar function.

\end{abstract}

    \maketitle

	\tableofcontents
	
	\newpage

\section{Introduction and outline} \label{sec:intro}

	Fractional calculus and the theory of special functions have managed to attract an increasing attention from the mathematical and physical community, especially in the last fifty years. In particular, the strict connection between these two research topics has been acting as the driving force for the latest developments in the literature on these subjects. The aim of this work is to build on these premises, collecting and summarizing the main works and results on the emerging Prabhakar's approach to fractional calculus. In this spirit, in the following we provide the reader, who already has a general understanding of the basics of fractional calculus, with a comprehensive practical guide to the main aspects of Prabhakar calculus. In other words, the target audience for this review consists of researchers with an elementary background in fractional calculus who need a quick and hands-on training in order to start working on this compelling topic.

	The review is therefore organized as follows: First, in Section \ref{sec:history} we summarize the key historical events that have led to the discovery and modern developments on the Mittag--Leffler function, its generalizations, and the extension of fractional calculus based on these exotic objects. In Section \ref{sec:motivations} we provide some physical motivations for the need of this new modeling scheme. In Section \ref{sec:math} we then recap the main mathematical features and properties of the Prabhakar function. In Section \ref{sec:pracalculus} we use the discussion in Section \ref{sec:math} as a springboard for the rigorous description of the theory of fractional calculus based on the Prabhakar function. In Section \ref{sec:physics} we profit of the preliminary mathematical sections in order to frame some anomalous physical phenomena, mentioned in Section \ref{sec:motivations}, within this generalized theory of calculus. In Section \ref{sec:probability} we analyze the implications that these operators bring onto the theory of probability, with particular attention for renewal processes. In Section \ref{sec:numerics} we provide a bird's-eye view on the numerical methods for evaluating the Prabhakar function. Finally, in Section \ref{sec:conclusions} we provide some concluding remarks and an overview on some open problems.
	
\vspace*{-3pt} 
	\section{History of the Mittag--Leffler function and its generalizations} \label{sec:history}

\setcounter{section}{2} \setcounter{equation}{0} 

	In the beginning of the twentieth century the Swedish mathematician Magnus G\"osta Mittag--Leffler, while working on summation methods for divergent series, defined a new special function which soon became known as the {\em Mittag--Leffler} (ML) {\em function}. Originally proposed in \cite{Mittag-Leffler1902,MittagLeffler1904}, the ML function was then analyzed in a systematic fashion by A. Buhl in \cite{Buhl1925} around twenty years from its original formulation. Most notably, this function turns out to be a one-parameter version of a more general (two-parameter) function promoted by A. Wiman in \cite{Wiman1905}, just a few years after the seminal works of M. G. Mittag--Leffler. Later on, in 1948, Pollard \cite{Pollard1948} investigated the complete monotonicity of the (one-parameter) ML function, while other efforts aimed at a deeper understanding of the main features of this function and its two-parameter generalization were made by P. Humbert, {\em et al.} in 1953 \cite{Humbert1953,HumbertDelerue1953,HumbertAgarwal1953}. Then, the first extensive collection of results on ML functions, with one and two parameters, appeared in 1955 in the third volume of the Bateman Manuscript Project \cite{ErdelyiMagnusOberhettingerTricomi3}. This project, aimed at redacting an encyclopedia on special functions, was started by Harry Bateman. However, the English mathematician passed away before the completion of this endeavor, which was ultimately brought to its end by A. Erd\'{e}lyi, W. Magnus, F. Oberhettinger, and F. G. Tricomi. At that time the role of the ML function and its generalizations was yet to be recognized in the literature. For this reason these functions only appear in the Bateman series as a miscellaneous topic. The current notoriety of the ML function and its generalization can however be traced back to the 1930s, specifically to the work of E. Hille and J. D. Tamarkin \cite{HilleTamarkin1930}, when the connection between these objects and the solution of weakly singular linear Volterra integral equations with power law kernels was drawn for the first time. Along these lines, in 1954 J. H. Barrett employed the ML function to express the solution of some linear differential equations of {\em non-integer order} \cite{Barrett1954}, thus identifying the fundamental role of this function in {\em fractional calculus}.

	The increasing interest in the theory of differential equations of non-integer order thus motivated the mathematical community to investigate at a deeper level the properties and features of the ML function and its generalization, as it is clearly highlighted by several Era-defining monographs
(see, \eg \cite{Djrbashian1993,KilbasSrivastavaTrujillo2006,Mainardi2010,MillerRoss1993,Podlubny1999}) and fundamental works (see, \eg \cite{MainardiGorenflo2000,KilbasTrujillo2001,Krageloh2003,HauboldMathaiSaxena2011}). Nonetheless, the first and only consistent and comprehensive monograph \cite{GorenfloKilbasMainardiRogosin2014} dedicated solely to this important topic was published in 2014 as the result of a joint endeavor of R. Gorenflo, A. A. Kilbas, F. Mainardi, and S. Rogosin.

	Aside from the (mostly) mathematical interest in the ML function, in the early days of fractional calculus, it is worth remarking that this subject had also a significant impact in the physics literature. This is particularly highlighted by the early works of K. S. Cole on the electrical conductance of biological systems \cite{Cole1933electric} (see also \cite{Davis1936theory}) and by the fundamental contributions of B. Gross to the study of dielectric and mechanical relaxation, see \eg \cite{Gross1947creep1,Gross1947creep2,Gross1947creep3}. Afterwards, M. Caputo and F. Mainardi gave their crucial contributions to fractional viscoelasticity \cite{CaputoMainardi1971-1,CaputoMainardi1971-2} by showing that the ML function always appears in the material functions of a system described by a fractional-order constitutive equation. It is also worth noting that \cite{CaputoMainardi1971-1} contains the very first plot of the ML function appeared in the literature. These are just a few early studies on the implications of the ML function in applied sciences that served as the backbone for all the subsequent developments in fractional calculus. Indeed, in light of the striking relevance of this function for the mathematical modeling of several phenomena, in 2007 the ML function was named the {\em Queen function of fractional calculus} \cite{MainardiGorenflo2007} by F. Mainardi and R. Gorenflo. Since we do not have the pretense  of providing a full and complete picture of the history and implications of this topic, we refer the interested reader to \cite{GorenfloKilbasMainardiRogosin2014} by R. Gorenflo, {\em et al.}

	The clear impact of the ML function and its two-parameter version on both the physical and mathematical community has then stimulated the quest for the generalization of these objects beyond two parameters. Among the many attempts we wish to commend the early contributions of V. Kiryakova 
and Y. Luchko 
on the so-called multi-index Mittag--Leffler functions \cite{Kiryakova1999}, \cite{Luchko1999}, \cite{Kiryakova2000}, \cite{Kiryakova2010}. 
If not for a few exceptions, among which we find the works of V. Kiryakova mentioned above, most of the proposed generalizations were just curiosity-driven and in general do not carry any relevance for physics or real-world applications. Contrarily, another special praise is due to the three-parameter extension of the ML function introduced by the Indian mathematician Tilak Raj Prabhakar in \cite{Prabhakar1971}, dating back to 1971. Although this function was originally analyzed within the framework of weakly singular Volterra integral equations from a pure mathematical perspective, it turns out that it plays an important role in the description of anomalous dielectric relaxation. An example of the effectiveness of the {\em Prabhakar function} in this context is provided by the Havriliak--Negami empirical law \cite{HavriliakNegami1967,HavriliakNegami1969a,HavriliakNegami1969b}. Notably, despite the short time-span between the original proposal of T. R. Prabhakar and the empirical studies of S. Havriliak and S. Negami, the connection between these works was drawn not longer than a decade ago, with the works of A. Hanyga and M. Seredy\'{n}ska \cite{HanygaSeredynska2008} together with the study of E. Capelas De Oliveira, {\em et al.} \cite{CapelasMainardiVaz2011}. Nonetheless, it is worth mentioning that several authors had previously investigated the time-domain representation of the Havriliak--Negami model without highlighting (explicitly) its relation with the Prabhakar function (see, \eg \cite{NigmatullinRyabov1997,Hilfer2002a,Hilfer2002b,Uchaikin2003,WeronJurlewiczMagdziarz2005,NigmatullinOsokinSmith2003,Stanislavsky2007,
NovikovWojciechowskiKomkovaThiel2005}).

	The first contribution aimed at providing a fractional interpretation for weakly singular Volterra operators involving the Prabhakar function was first developed in 2002 by A. A. Kilbas, M. Saigo, and R. K. Saxena in \cite{KilbasSaigoSaxena2002}. Specifically, in this paper the notion of the Prabhakar fractional integral was first identified as an independent entity worth of investigation. Besides, in \cite{KilbasSaigoSaxena2004} one of the left-inverse operators of the Prabhakar fractional integral, which will later become known as Prabhakar fractional derivative (of the Riemann-Liouville type), was found. Over a decade later, in 2013, M. D'Ovidio and F. Polito \cite{DOvidioPolito2018} provided a regularization of the Prabhakar fractional derivative and named it after the Indian mathematician, thus originating the current terminology on the subject. These preliminary analyses  then inspired the seminal work by R. Garra, {\em et al.} \cite{GarraGorenfloPolitoTomovski2014} where the features of the Prabhakar derivatives and some of their applications are discussed in depth. This work, in particular, contributed to kick-start a series of studies concerning applications of these operators to different fields of research (see, \eg \cite{ChamatiTonchev2006,SaxenaMathaiHaubold2006,CamargoChiacchioCharnetCapelasDeOliveira2009,
CamargoCharnetCapelasDeOliveira2009,SandevTomovski2010,TomovskiHilferSrivastava2010}). A particular praise is due to the work of L. Beghin and E. Orsingher who first derived the connection between the Prabhakar and Wright functions \cite{BeghinOrsingher2010}. Furthermore, it is of paramount importance to mention the contribution of R. B. Paris \cite{Paris2010} in the study of the asymptotic behavior of the Fox--Wright functions, of which the Prabhakar function represents a special case, which has then led to specific results on the Prabhakar function \cite{GarraGarrappa2018,Paris2019}.

	Finally, we conclude this bird's-eye view on the history of this subject by mentioning some recent results carrying some more practical implications. First, the conditions for the complete monotonicity of the Prabhakar function have been extensively analyzed throughout the last decade in \cite{CapelasMainardiVaz2011,MainardiGarrappa2015_JCP,TomovskiPoganySrivastava2014} by several authors. Using these studies A. Giusti \cite{Giusti2019} has then laid down a connection between the so called Prabhakar fractional calculus and Kochubei's general fractional calculus \cite{Kochubei2011}. Lastly, it is worth recalling that a discussion of the numerical aspects of the Prabhakar function has been carried out by R. Garrappa in \cite{Garrappa2015_SIAM}, who also designed a Matlab code for computing the Prabhakar function \cite{Garrappa2014ml}. Note that, as of today, this is the sole code (freely available on the Mathwork website) for the computation of the Prabhakar function.

\vspace*{-3pt} 
\section{Anomalous physics: a cry for help} \label{sec:motivations}

\setcounter{section}{3} \setcounter{equation}{0} 

	In this section we provide some physical examples showing the need for an extension of ordinary calculus based on the Prabhakar function. In this regard, the most telling framework is provided by the so called \emph{anomalous phenomena} emerging in several physical settings.
	
	Let us start off with the classical theory of dielectrics \cite{Jackson}. A dielectric material is a poor conductor of electricity whose main feature is that its internal electric charges are not allowed to flow freely under the effect of an external electric field. On the contrary, they will mostly tend to shift a little from their equilibrium position. This will result in a very low electric current, since some charges will still be able to move through the material, and a \emph{dielectric polarization} of the system. Experimentally it is observed that once the external stress is turned off a dielectric material will tend to lose its acquired polarization over time. This time lag, known in physics as \emph{dielectric relaxation}, is due to the fact that the atomic and/or molecular polarization in response to a changing electric field is not an instantaneous process (see \eg \cite{Jonscher1999}). A key quantity in the study of electromagnetism in matter is given by the \emph{polarization density} $\bm{P}$, or simply \emph{polarization}, which is defined as the average dipole moment of the material per unit volume and it accounts for the contribution of the ``bound charges'' that are responsible for the dielectric polarization. This picture allows us to split the electric charge density appearing in Gauss' law into two contributions, namely the one given by the charges which are free to move through the material and a second piece due to bound charges. This suggests the introduction of a new entity
\vskip -12pt
$$ \bm{D} = \epsilon_0 \, \bm{E} + \bm{P}, $$
known as the \emph{electric displacement}, where $\bm{E}$ is the electric field and $\epsilon_0$ denotes the \emph{vacuum permittivity}, whose divergence depends solely on the free charge.

	For a linear dielectric one has that the polarization density is proportional to the electric field, namely
\vskip -10pt
	\be \label{eq-linearPE}
	\bm{P} = \epsilon_0 \, \chi \, \bm{E} \, ,
	\ee
\vskip -4pt \noindent
	which leads to
\vskip -12pt
	\be
	\bm{D} = \epsilon \, \bm{E} = \epsilon_0 \, \epsilon_{\rm r} \, \bm{E} =
	\epsilon_0 \, (1 + \chi) \, \bm{E} \, ,
	\ee
	with $\chi$ the \emph{electric susceptibility} of the material, $\epsilon_{\rm r} = 1 + \chi$ the \emph{relative permittivity}, and $\epsilon = \epsilon_0 \, \epsilon_{\rm r}$ the \emph{permittivity} of the system. On the other hand, it is well known that not all dielectrics display this very simple behavior, hence the polarization density can be, in general, a complicated function of the electric field. Nonetheless, if one wishes to keep the relation between $\bm{P}$ and $\bm{E}$ linear while introducing a non-instantaneous response to a time-varying electric field, then a straightforward generalization of \eqref{eq-linearPE} is obtained via a Fourier convolution. In detail, given two sufficiently regular functions $f(t)$ and $g(t)$ the convolution of these two functions, in the Fourier sense, is defined as
\vskip -12pt
	\be
	(f \ast g) (t) := \int _{-\infty} ^{\infty} f(t-\tau) \, g(\tau) \, \d \tau \, .
	\ee
	
Hence, if we now denote by $\chi (t)$ a time dependent electric susceptibility such that $\chi (t) = 0$ for $t<0$, then \eqref{eq-linearPE} can be generalized as
\vskip -10pt
	\be \label{PconvE}
	\bm{P} (t) = \epsilon_0 \left( \chi \ast \bm{E} \right) (t) \equiv
	\epsilon_0 \int _{-\infty} ^t \chi (t -\tau) \, \bm{E} (\tau) \, \d \tau \, ,
	\ee
\vskip -4pt \noindent
	or, in the frequency domain
\vskip -11pt
	\be
	\wh{\bm{P}} (\omega) = \epsilon_0 \, \wh{\chi} (\omega) \, \wh{\bm{E}} (\omega) \, ,
	\ee
\vskip -3pt \noindent
	where the hat denotes the Fourier transform, \ie
	$$
	\wh{f} (\omega) \equiv \mathcal{F}\big[ f(t) \, ; \, \omega \big] :=
	\int _{-\infty} ^{\infty} \exp \left( - i \, \omega \, t \right) \, f(t) \, \d t \, .
	$$
		
	Eq. \eqref{PconvE} is particularly instructive since it highlights the role of $\chi (t)$ as a measure of the response of the material to a sudden variation in the electric field $\bm{E}(t)$. In other words, the time dependent electric susceptibility acts as a ``response function'' of the dielectric which, in this scenario, is pictured as a passive linear system (see \eg \cite{GarrappaMainardiMaione2016,Vaughan1979}, and references therein).  More precisely, the \emph{response function} of the dielectric can be defined, in the frequency domain, as
\vskip -12pt
	\be \label{transferphi}
	\wh{\phi} (\omega) :=
	\frac{\wh{\chi}(\omega)}{\wh{\chi}(0)} \, .
	\ee

Now, since we are considering a causal linear system and we are only interested in the history of the system from a certain time $t=0$ (at which we have applied the external perturbation) up to the time $t>0$, Eq. \eqref{transferphi} can be recast as \cite{GarrappaMainardiMaione2016}
\vskip -12pt
	\be
	\wt{\phi} (s) = \wh{\phi} (\omega) \big|_{s=i \omega} \, ,
	\ee
	where the tilde denotes the Laplace transform, \ie
\vskip -10pt
	$$
	\wt{f} (s) \equiv \mathcal{L}\big[ f(t) \, ; \, s \big] :=
	\int _{0} ^{\infty} \exp \left( - s \, t \right) \, f(t) \, \d t \, .
	$$
	
	Alongside with the response function one can also introduce another quantity, known as  \emph{relaxation function}, which is defined as
	\be \label{relaxation-function}
	\Psi (t) := 1 - \int _0 ^t \phi (\tau) \, \d \tau \, , \qquad t \geq 0 \, ,
	\ee
	from which one can also infer that $\phi (t) = - \Psi ' (t)$, where the prime denotes the derivative with respect to time.

	The standard lore on dielectrics is described by the \emph{Debye model} \cite{Debye-Seminal,Debye} according to which an ideal dielectric is schematically represented as a series of non-interacting dipoles. Such a model then leads to a (normalized) complex susceptibility given by
\vskip -11pt
	\be
	\wh{\chi}  _{\rm D}(\omega) = \frac{1}{1 + \iu \, \omega \, \tau _{\rm D}} \, ,
	\ee
\vskip -4pt \noindent
	which yields \cite{GarrappaMainardiMaione2016, Debye}
\vskip -10pt
	\be
	\phi _{\rm D} (t) = \frac{\exp (- t / \tau_{\rm D})}{\tau_{\rm D}} \, , \quad
	\Psi _{\rm D} (t) = \exp( - t / \tau_{\rm D}) \, ,
	\ee
\vskip -3pt \noindent
	with $\tau_{\rm D}$ denoting the relaxation time.
	
	While this very simple and powerful model provides a satisfactory description of the features of a wide class of standard dielectrics, there are several experimental evidences \cite{Jonscher1999,exp1,exp2} that show large deviations from the exponential relaxation for many materials. An example of these \emph{anomalous dielectrics} is provided by the \emph{Cole--Cole relaxation} \cite{CC1,CC2} displayed by certain biological tissues \cite{CC-exp1,CC-exp2}. The key feature of these systems is a fractional-power law fall-off shown by the complex susceptibility, at high frequencies, which leads to a stretched relaxation over a wider range of frequencies. In detail, the Cole--Cole relaxation is modeled in terms of a complex susceptibility that reads
\vskip -11pt
	\be
	\wh{\chi} _{\rm CC} (\omega) = \frac{1}{1 + (\iu \, \omega \, \tau _{\star})^\alpha} \, ,
	\quad 0 < \alpha \leq ,1
	\ee
\vskip -4pt \noindent
that yields
\vskip -14pt
	\be \label{phi-CC}
	\phi _{\rm CC} (t) = \frac{1}{\tau_{\star}} \left(\frac{t}{\tau_{\star}} \right) ^{\alpha - 1}
	E_{\alpha, \alpha} \left[ - \left(\frac{t}{\tau_{\star}} \right) ^{\alpha} \right]
	\ee
\vskip -4pt \noindent
	and
\vskip -14pt
	\be \label{psi-CC}
	\Psi _{\rm CC} (t) =
	E_{\alpha, 1} \left[ - \left(\frac{t}{\tau_{\star}} \right) ^{\alpha} \right] \, ,
	\ee
where $\tau _\star$ represents a typical relaxation timescale for the system.

	Note that, although the Cole--Cole model breaks the standard exponential relaxation replacing it with fractional-power law tails, it can still be framed within the conventional picture of fractional calculus based on Caputo derivatives. This is, however, not the case for several other empirical laws for anomalous dielectrics such as the \emph{Davidson--Cole} \cite{DC} and \emph{Havriliak--Negami} \cite{HN} models, described by \cite{GarrappaMainardiMaione2016}
\vskip -12pt
	\be \label{DC}
	\wh{\chi} _{\rm DC} (\omega) = \frac{1}{(1 + \iu \, \omega \, \tau _{\star})^\gamma} \, ,
	\quad 0 < \gamma \leq 1
	\ee
\vskip -4pt \noindent
and \cite{MainardiGarrappa2015_JCP, GarrappaMainardiMaione2016, HH1994}	
\vskip-13
pt		
	\be \label{HN}
	\wh{\chi} _{\rm HN} (\omega) = \frac{1}{[1 + (\iu \, \omega \, \tau _{\star})^\alpha]^\gamma} 		\, ,
	\quad 0 < \alpha \leq 1 \, , \,\,\, 0 < \alpha \gamma \leq 1
	\ee			
respectively, for which the dielectric decay is inherently related to the Prabhakar function \cite{GarrappaMainardiMaione2016,MainardiGarrappa2015_JCP,GarrappaHN2016,CapelasMainardiVaz2011,
GorskaJPA2018,GorskaPLA2019}, as we shall discuss in detail in Section \ref{sec:physics}.		
			
	Although anomalous dielectric relaxation is probably the most apparent exemplification of the need for a fractional theory of calculus based on the Prabhakar function, there exist also several other different physical systems in which this formalism naturally emerges. For instance, this scenario comes up when dealing with certain kind of anomalous diffusion processes \cite{Sandev-EPL,Sandev-Diffusion1,Sandev-Diffusion2,Sandev-Diffusion3,Sandev-Diffusion4,Sandev-Diffusion5,Sandev-Diffusion6,Sandev-Diffusion7}, and when employing the formal duality between dielectric materials and viscoelastic systems \cite{GC-CNSNS2018,CGV-Mathematics2018,Giusti2018_NODY}. More on these physical cases can be found in Section \ref{sec:physics}.

\section{Mathematical preliminaries on the Prabhakar function} \label{sec:math}

\setcounter{section}{4} \setcounter{equation}{0} 

The Prabhakar function (also known as the three-parameter Mittag-Leffler function, see \eg \cite{GorenfloKilbasMainardiRogosin2014}), introduced by T. R. Prabhakar in 1971 \cite{Prabhakar1971}, is defined as
\begin{equation}\label{eq:PrabhakarFunction}
		E_{\alpha,\beta}^{\gamma}(z) = \sum_{k=0}^{\infty} \frac{(\gamma)_k \, z^{k}}{k! \Gamma(\alpha k + \beta)}
		, \quad
		\alpha, \beta, \gamma \in \mathbb{C}
		, \quad \real(\alpha) > 0 , \,\,\,  z \in \C \, ,
\end{equation}
where $(\gamma)_k \equiv \Gamma (\gamma + k) / \Gamma (\gamma)$ is the Pochhammer symbol and $\Gamma(\cdot)$ denotes the Euler gamma function. Note that $E_{\alpha,\beta}^{\gamma}(z)$ is an {\em entire function} of order $\rho = 1/\real(\alpha)$ and type $\sigma=1$ \cite{GorenfloKilbasMainardiRogosin2014}.
	
Although the three parameters $\alpha$, $\beta$ and $\gamma$ are allowed to assume values in $\C$, and several of the properties we are going to present hold with complex parameters, in this work we will restrict our attention just to real parameters. The reason for this choice is that real parameters turn out to be of particular interest in physics, as we shall discuss later on. Therefore, for the remainder of this work we will implicitly assume
\[
	\alpha, \beta, \gamma \in \R, \quad \alpha > 0 .
\]
	
	In view of the increasing relevance of the Prabhakar function in current literature, in the following we summarize the most relevant results and open problems on this compelling research topic.
	
\subsection{Main properties and relations with other functions} 

	The Prabhakar function has important connections with the standard ML function, its two-parameter generalization, and other special functions. Since for $k\ge 1$ one has
\vskip-10pt
$$
(\gamma)_k = \frac{\Gamma(\gamma+k)}{\Gamma(\gamma)} = \gamma (\gamma+1)\cdots(\gamma+k-1) \, ,
$$
then for $\gamma=0$ we clearly have
\vskip -10pt
$$
	E_{\alpha,\beta}^{0}(z) = \frac{1}{\Gamma(\beta)} \, ,
$$
whereas, for $\gamma=1$ we recover the widely known two-parameter ML function, \ie
\vskip -12pt
\[
	E_{\alpha,\beta}^{1}(z) = E_{\alpha,\beta}(z) = \sum_{k=0}^{\infty} \frac{z^k}{\Gamma(\alpha k + \beta)} \, .
\]

Hence, if one sets $\alpha = \beta = \gamma = 1$ the classical exponential function
$
E_{1,1}^{1}(z) = \eu^{z}
$ is retrieved.

	Being the Prabhakar function distinguished from the two-parameter ML function just for the presence of a third parameter $\gamma$, one is compelled to focus on properties related to this additional parameter.

A first formula, already present in the original work by Prabhakar \cite{Prabhakar1971}, enables one to operate the reduction of the third parameter, \ie
\begin{equation}\label{ReductionFormula1}
	E_{\alpha,\beta}^{\gamma+1}(z) =
	\frac{E_{\alpha,\beta-1}^{\gamma}(z) + (1-\beta+\alpha\gamma) E_{\alpha,\beta}^{\gamma}(z)}{\alpha \gamma} ,
\end{equation}
and a further reduction formula was later derived in \cite{GarraGarrappa2018} and reads
\begin{equation}\label{ReductionFormula2}
	E_{\alpha,\beta}^{\gamma+1}(z) =
	\frac{E_{\alpha,\beta-\alpha-1}^{\gamma}(z) + (1-\beta+\alpha) E_{\alpha,\beta-\alpha}^{\gamma}(z)}{\alpha \gamma z} ,
	\quad z \not= 0 .
\end{equation}

Both these formulas can be used to lower the value of the third (additional) parameter, nonetheless these results can probably be better appreciated when $\gamma$ is an integer. Specifically, let $\gamma = k \in \N$, then
\vskip -10pt
\begin{equation}\label{eq:PrabhakarSummationFormula}
	E_{\alpha,\beta}^{k+1}(z) = \frac{1}{\alpha^k k!} \sum_{j=0}^{k} d_{j}^{(k)} E_{\alpha,\beta-j}(z)
	\, ,
\end{equation}
\vskip -4pt \noindent
and
\vskip -13pt
\begin{equation}\label{eq:DzhrbashyanSummationFormula}
	E_{\alpha,\beta}^{k+1}(z) = \frac{1}{\alpha^k z^{k}k!} \sum_{j=0}^{k} d_{j}^{(k)} E_{\alpha, \beta - \alpha k - j}(z)
	\, , \quad z \neq0 \, ,
\end{equation}
which means that one can explicitly express $E_{\alpha,\beta}^{k+1}$ as a combination of two-parameter ML functions. The coefficients $d_{j}^{(k)}$ appearing in both (\ref{eq:PrabhakarSummationFormula}) and (\ref{eq:DzhrbashyanSummationFormula}) are given by the recursive expression
\begin{equation}\label{eq:CoeffSummationFormula}
	d_{j}^{(k)} =
	\left\{ \begin{array}{ll}
		(1 - \beta + \alpha ) d_{0}^{(k-1)} & j = 0 ,\\
		d_{j-1}^{(k-1)} + (1 - \beta + \alpha + j) d_{j}^{(k-1)} \quad & j=1,\dots,k-1 ,\\
		1 & j = k .\
	\end{array}\right.
\end{equation}
originally derived by R. Garrappa and M. Popolizio in \cite{GarrappaPopolizio2018}.
	
	The Prabhakar function is also directly related to the Fox--Wright functions \cite{Fox1928,Wright1935}. Specifically, the latter are defined as
\begin{equation}\label{FoxWright}
		_{p}\Psi_{q}(z) \equiv \,
		_{p}\Psi_{q} \left[ \begin{array}{c} (a_1, \rho_1),\dots,(a_p, \rho_p) \\
		(b_1,\sigma_1),\dots,(b_q, \sigma_q) \end{array} ; z \right]
		= \sum_{k=0}^{\infty}
			\frac
			{\prod_{r=1}^{p} \Gamma(a_{r} + \rho_{r} k) }
			{\prod_{s=1}^{q} \Gamma(b_{s} + \sigma_{s} k) }
			\frac{z^{k}}{k!} \, ,
\end{equation}
where $p$ and $q$ are integers, $\rho_r, a_r, \sigma_r,b_r$ are real or complex parameters such that $\rho_r k + a_r \neq 0, -1, -2, \ldots$. Then, one can easily verify that $E_{\alpha,\beta}^{\gamma}(z)$ is proportional to $_{1}\Psi_{1}(z)$ since
\vskip -10pt
	\begin{equation}\label{eq:PrabhakarFoxWright}
		E_{\alpha,\beta}^{\gamma}(z)
		= \frac{1}{\Gamma(\gamma)} \, _{1}\Psi_{1} \left[ \begin{array}{c} (\gamma,1) \\ (\beta,\alpha) \end{array} ; z 			\right] .
	\end{equation}

Moreover, if we recall the relation between the Fox--Wright functions and the Fox $H$--function \cite{mathai2009h}, \ie
\begin{equation}\label{eq:PrabhakarFoxH}
\begin{split}
{}_{p}\Psi_{q} &\Big[ \begin{array}{c} (a_1, \rho_1),\dots,(a_p, \rho_p) \\
		(b_1,\sigma_1),\dots,(b_q, \sigma_q) \end{array} ; z \Big]\\
 & = H^{1,p} _{p, q+1}
 \left[ - z \left|
		\begin{array}{c}
		(1 - a_1, \rho_1), \, (1 - a_2, \rho_2), \dots,(1 - a_p, \rho_p) \\
		(0,1), \, (1- b_1,\sigma_1),\dots,(1 - b_q, \sigma_q) \end{array}
		\right. \right] \, ,
\end{split}
\end{equation}
then \eqref{eq:PrabhakarFoxWright} can be rewritten in terms of the Fox $H$-function as (see \cite{SaxenaGianni2011})
\begin{equation}
E_{\alpha,\beta}^{\gamma}(z) =
\frac{1}{\Gamma(\gamma)}
H^{1,1} _{1,2}
\left[ - z \left|
		\begin{array}{c}
		(1 - \gamma, 1)  \\
		(0,1), \, (1- \beta,\alpha)
		\end{array}
		\right. \right] \, .
\end{equation}

Finally, it is worth mentioning that series in Prabhakar functions have been studied by J. Paneva--Konovska in \cite{paneva2013multi,paneva2014convergence,paneva2017overconvergence}.

\subsection{Derivatives and integrals} 

A term-by-term derivation of the series expansion in \eqref{eq:PrabhakarFunction} gives
\begin{equation*}
\begin{split}
	\frac{\du^m}{\du z^m} E_{\alpha,\beta}^{\gamma}(z) &= \gamma(\gamma+1)\cdots(\gamma+m-1) E_{\alpha,m\alpha+\beta}^{\gamma+m}(z)\\
	&= \frac{\Gamma(\gamma+m)}{\Gamma(\gamma)} E_{\alpha,m\alpha+\beta}^{\gamma+m}(z) \, ,
\end{split}
\end{equation*}
with $m \in \N$, that reduces to the known result for the repeated derivative of the two-parameter ML function for $\gamma = 1$, \ie
\begin{equation}\label{eq:DerMLPrabhakar}
	\frac{\du^m}{\du z^m} E_{\alpha,\beta}(z) = m! \, E_{\alpha,m\alpha+\beta}^{m+1}(z) \, .
\end{equation}

This last expression turns out to be particularly useful when combined with (\ref{eq:PrabhakarSummationFormula}) or (\ref{eq:DzhrbashyanSummationFormula}), since it allows to express derivatives of the two-parameter ML function as combination of different instances of the same function, a result widely exploited in \cite{GarrappaPopolizio2018} for numerical purposes.

	Similarly, {\em Dzhrbashyan's formula} \cite{Djrbashian1993, GorenfloKilbasMainardiRogosin2014} for the derivative of the two-parameter ML function can be generalized to the Prabhakar case as (see \cite{GarraGarrappa2018})
\vskip -12pt
\begin{equation}\label{DerivPrabhakar2}
	\frac{\du}{\du z} E_{\alpha,\beta}^{\gamma}(z) = \frac{E_{\alpha,\beta-1}^{\gamma}(z) + (1-\beta)E_{\alpha,\beta}^{\gamma}(z)}{\alpha z} \, , \quad z\neq 0 \, .
\end{equation}

	Other important results concerning integrals and derivatives involving the Prabhakar function appear to be strictly related to the product $t^{\beta-1} E_{\alpha,\beta}^{\gamma}(t^{\alpha} z)$. In details, it is easy to see that the $m$-th derivative and the integral of this particular combination read
\vskip -10pt
\[
	\frac{\du^{m}}{\du t^{m}} t^{\beta-1} E_{\alpha,\beta}^{\gamma}(t^{\alpha} z) = t^{\beta-m-1} E_{\alpha,\beta-m}^{\gamma}(t^{\alpha} z) \, , \quad m \in \N
\]
\vskip -4pt \noindent
and
\vskip -10pt
\[
	\int_{0}^{t} \tau^{\beta-1} E_{\alpha,\beta}^{\gamma}(\tau^{\alpha} z) \du \tau = t^{\beta} E_{\alpha,\beta+1}^{\gamma}(t^{\alpha} z) \, ,
\]
respectively. Further, the $m$-th repeated integration of $t^{\beta-1} E_{\alpha,\beta}^{\gamma}(t^{\alpha} z)$ gives
\vskip -10pt
\be
\notag J_0 ^m \big[ t^{\beta-1} E_{\alpha,\beta}^{\gamma}(t^{\alpha} z) \big] &\equiv&
\int_0^t \du \tau_1 \int_0^{\tau_1}\du \tau_2 \cdots \int_0^{\tau_{m-1}}\du \tau_m \, \tau_m^{\beta-1} E_{\alpha,\beta}^{\gamma}(\tau_m^{\alpha} z) \\
&=& \notag
\frac{1}{(m-1)!} \int_{0}^{t} (t-\tau)^{m-1}
\tau^{\beta-1} E_{\alpha,\beta}^{\gamma}(\tau^{\alpha} z) \d \tau \\
&=&
t^{\beta+m-1} E_{\alpha,\beta+m}^{\gamma}(t^{\alpha} z) \, .
\ee

	We conclude this section by collecting some formulas for fractional integrals and derivatives of the Prabhakar function. To this aim, we first need to recall that for a function $f \in L^1[t_0,T]$ the {\em Riemann-Liouville (RL) fractional integral} of order $\alpha>0$ is defined as \cite{Diethelm2010,KilbasSrivastavaTrujillo2006,Mainardi2010,MillerRoss1993,Samko-book}
\begin{equation}\label{eq:RL_Integral}
	J^{\alpha}_{t_0} f(t) = \frac{1}{\Gamma(\alpha)} \int_{t_0}^{t} (t-u)^{\alpha-1} f(u) \, \du u \, .
\end{equation}

One of its left-inverse operators, known as the {\em Riemann-Liouville fractional derivative}, then reads
\vskip -12pt
\begin{equation}\label{eq:RL_Derivative}
	\DR^{\alpha}_{t_0} f(t)
	\coloneqq D^{m} J^{m-\alpha}_{t_0} f(t)
	=  \frac{1}{\Gamma(m-\alpha)} \frac{\du^m}{\du t^m} \int_{t_0}^{t} (t-u)^{m-\alpha-1} f(u) \du u
\end{equation}
with $m = \left\lceil \alpha \right\rceil$ (the smallest integer greater than $\alpha$). A regularized version of the RL derivative, which still acts as left-inverse of the RL integral, is given by the {\em Caputo fractional derivative}. The latter is obtained by exchanging the operations of integration and differentiation in (\ref{eq:RL_Derivative}) and applies to functions $f \in AC^m[t_0,T]$ (the set of functions with absolutely continuous derivatives of order $m-1$), \ie
\vskip -12pt
\begin{equation}\label{eq:C_Derivative}
	\DC^{\alpha}_{t_0} f(t)
	\coloneqq J^{m-\alpha}_{t_0} D^{m}  f(t)
	=  \frac{1}{\Gamma(m-\alpha)} \int_{t_0}^{t} (t-u)^{m-\alpha-1} f^{(m)}(u) \, \du u \, .
\end{equation}

All that being said, it is straightforward to verify that (see \eg \cite{KilbasSaigoSaxena2004,GorenfloKilbasMainardiRogosin2014})
\[
	\begin{aligned}
	J^{\rho}_{0} \left[ t^{\beta-1} E_{\alpha,\beta}^{\gamma}\bigl( t^{\alpha} \lambda \bigr) \right]
	&=
	t^{\beta+\rho-1} E_{\alpha,\beta+\rho}^{\gamma}\bigl( t^{\alpha} \lambda\bigr)  \\
	\DR^{\rho}_{0} \left[ t^{\beta-1} E_{\alpha,\beta}^{\gamma}\bigl( t^{\alpha} \lambda \bigr) \right]
	&=
	t^{\beta-\rho-1} E_{\alpha,\beta-\rho}^{\gamma} \bigl(t^{\alpha} \lambda \bigr) \\
	\DC^{\rho}_{0} \left[ t^{\beta-1} E_{\alpha,\beta}^{\gamma}\bigl( t^{\alpha} \lambda \bigr) \right]
	&=
	t^{\beta-\rho-1} E_{\alpha,\beta-\rho}^{\gamma}\bigl( t^{\alpha} \lambda \bigr) \, ,
	\quad \beta > \left\lceil \rho \right\rceil \, , \\
	\end{aligned}
\]
\vskip -3pt \noindent
with $\rho > 0$.

\subsection{Integral transforms} 

	Let us consider, again, the product $t^{\beta-1} E_{\alpha,\beta}^{\gamma}(t^{\alpha} \lambda)$. Then it is not hard to show that (see \eg \cite{KilbasSaigoSaxena2002})
\vskip -11pt
\begin{equation}\label{ML3_LT}
	{\mathcal L} \Bigl[ t^{\beta-1} E_{\alpha,\beta}^{\gamma}(t^{\alpha} z) \, ; \, s \Bigr]
= \frac{s^{\alpha\gamma - \beta}}{\bigl(s^{\alpha} - z\bigr)^{\gamma}}
	, \quad \real(s) > 0 \ \text{and} \ |s|>|z|^{\frac{1}{\alpha}} \, .
\end{equation}
This result is particularly useful since it plays an important role in both physical applications and numerical computations of the Prabhakar function. Alternatively, proceeding in a more general fashion one finds (see \cite{SrivastavaTomovski2009,SandevTomovski2010})
\vskip -13pt
\be
{\mathcal L} \Big[ t^{\rho-1} E_{\alpha,\beta}^{\gamma}(t^{\sigma} z) \, ; \, s \Big] = \frac{s^{-\rho}}{\Gamma(\gamma)} \, {}_2\Psi_{1} \Bigl[ \begin{array}{c} (\rho,\sigma), \, (\gamma,1) \\ (\beta,\alpha) \\ \end{array} ; \frac{z}{s^{\sigma}} \Bigr].
\ee

	Setting $t=1$ in \eqref{ML3_LT}, one finds an integral representation of the Prabhakar function, \ie
\vskip -13pt
\[
	E_{\alpha,\beta}^{\gamma}(z) = \frac{1}{2\pi \iu} \int_{\mathcal C} \eu^{s} \frac{s^{\alpha\gamma - \beta}}{\bigl(s^{\alpha} - z\bigr)^{\gamma}} \, \du s ,
\]	
\vskip -3pt \noindent
where ${\mathcal C}$ is understood as the {\em Bromwich contour}, for which all singularities lay on the left-hand side of a certain $\real(s) = \xi$ known as the {\em abscissa of convergence}. Similarly, a particular {\em spectral representation} of $t^{\beta-1} E_{\alpha,\beta}^{\gamma}(t^{\alpha} \lambda)$ has been computed in \cite{CapelasMainardiVaz2011,MainardiGarrappa2015_JCP,TomovskiPoganySrivastava2014}, and it reads
\vskip -11pt
\begin{equation}\label{eq:SpectralRepresentation}
	t^{\beta-1} E_{\alpha,\beta}^{\gamma}(-t^{\alpha}) = \int_{0}^{\infty} \eu^{-rt} K_{\alpha,\beta}^{\gamma}(r) \du r
\end{equation}
\vskip -4pt \noindent
with
\vskip -14pt
\be
	\notag K_{\alpha,\beta}^{\gamma}(r) &=&
	\frac{r^{\alpha\gamma - \beta}\sin\bigl( \gamma \theta_\alpha(r)  + (\beta - \alpha\gamma)\pi \bigr)}
	{\pi\bigl( r^{2\alpha} + 2 r^{\alpha} \cos(\alpha\pi) + 1 \bigr)^{\gamma/2}} \, ,\\
	\notag \theta_\alpha(r) &=& \arctan\left( \frac{r^\alpha \sin (\pi \alpha)}{r^\alpha \cos(\pi \alpha) +1} \right) \in [0,\pi] \, ,
\ee
\vskip -3pt \noindent
where a very specific {\em branch} of $\arctan (z)$, treated as a multi-valued function, has been chosen in order to have $\theta_\alpha(r) \in [0,\pi]$ (see \cite{MainardiGarrappa2015_JCP} for details). Further, this representation plays a key role in the determination of the range of parameters that yields a completely monotonic behavior of the Prabhakar function. Besides, this last representation is also particularly useful when studying estimates and bounds on the fractional Prabhakar integral, see \eg \cite{PolitoTomovski2016}.

In the light of the strict relation between Fox--Wright functions and the Prabhakar function, as shown in Eq. \eqref{eq:PrabhakarFoxWright}, one can easily derive the {\em Mellin--Barnes integral representation} of the three-parameter ML function, \ie
\be \label{eq:mellinbarnesrep}
	E_{\alpha,\beta}^{\gamma}(z) = \frac{1}{\Gamma(\gamma)} \frac{1}{2\pi \iu} \int_{{\mathcal C}} \frac{\Gamma(s)\Gamma(\gamma-s)}{\Gamma(\beta-\alpha s)} (-z)^{-s} \, \du s
	, \quad
	|\arg z| < \pi \, ,
\ee
where $\mathcal{C}$ is defined as above. Along this line, denoting by
\be
\mathcal{M} \Big[ f(t) \, ; \, s \Big] = \int _0 ^\infty f(t) \, t^{s-1} \, \d t \, ,
\ee
the {\em Mellin transform} of $f(t)$, then from \eqref{eq:mellinbarnesrep} one can immediately infer that
\vskip -10pt
\be
\mathcal{M} \Big[ E_{\alpha,\beta}^{\gamma} (-t) \, ; \, s \Big] =
 \frac{\Gamma (s) \Gamma (\gamma - s)}{\Gamma (\gamma) \Gamma (\beta - \alpha s)} \, .
\ee	
	
\subsection{Complete monotonicity} 
	
	A crucial role for the applications of special functions in relaxation models and probability theory is played by the notion of {\em complete monotonicity}. We recall that a function $f:(0,+\infty) \to \mathbb{R}$ is completely monotonic (CM) if $f$ has derivatives of all orders on $(0,+\infty)$ and $(-1)^{k}f^{(k)}(t) \ge 0$, for any $k\in \mathbb{N}\cup\{0\}$ and $t > 0$.
	
    The complete monotonicity of the Prabhakar function has been studied in several works \cite{CapelasMainardiVaz2011, MainardiGarrappa2015_JCP,TomovskiPoganySrivastava2014} and the current most general result states that the product $t^{\beta-1}E_{\alpha,\beta}^{\gamma}(-t^{\alpha})$ is CM if the three parameters satisfy the conditions
    \vskip -10pt
\begin{equation}\label{eq:CM_Range}
	0 < \alpha \le 1 , \quad 0 < \alpha \gamma \le \beta \le 1
\end{equation}
which are obtained by combining the representation in (\ref{eq:SpectralRepresentation}) and the {\em Bernstein theorem}. The latter, in particular, states that a function is CM if its spectral function is non-negative. Thus, the restrictions in \eqref{eq:CM_Range} are obtained by inspecting the conditions according to which $K_{\alpha,\beta}^{\gamma}(r) \ge 0$ for all $r\ge0$. This result is particularly relevant since $t^{\beta-1}E_{\alpha,\beta}^{\gamma}(-t^{\alpha})$ comes up very often in the study of anomalous dielectrics, as we shall discuss in more details in Section \ref{sec:physics}.

Finally, in \cite{GorskaHorzelaLattanziPogani2018} it was shown that $E_{\alpha,\beta}^{\gamma}(-t)$ is CM if $0<\alpha<1$ and $0 < \alpha \gamma \le \beta$, namely without the restriction $\beta \le 1$. 
For similar results on another Mittag-Leffler type function with 3 parameters, the so-called ``Le Roy type" function, see \cite{Gorska-etc-FCAA2019}.	
		
\subsection{Asymptotic behavior}\label{SS:AsymptoticBehavior}	
	
The asymptotic behavior of the Prabhakar function, in the whole complex plane, is a non-trivial topic due to the dependence of the coefficients in the asymptotic expansions on the three parameters $\alpha$, $\beta$ and $\gamma$ (see \eg \cite{GarraGarrappa2018} and references therein). On the other hand, rewriting the Prabhakar function as a Fox--Wright function, as in (\ref{eq:PrabhakarFoxWright}), allows one to take profit of the well-established results on the asymptotic behavior of this more general class of functions, see \eg \cite{Braaksma1964,Paris2010,Paris2019,Wright1935,Wright1940,Wrigth1940_PTRSL}.

\begin{figure}[ht]
\centering
\includegraphics[width=0.70\textwidth]{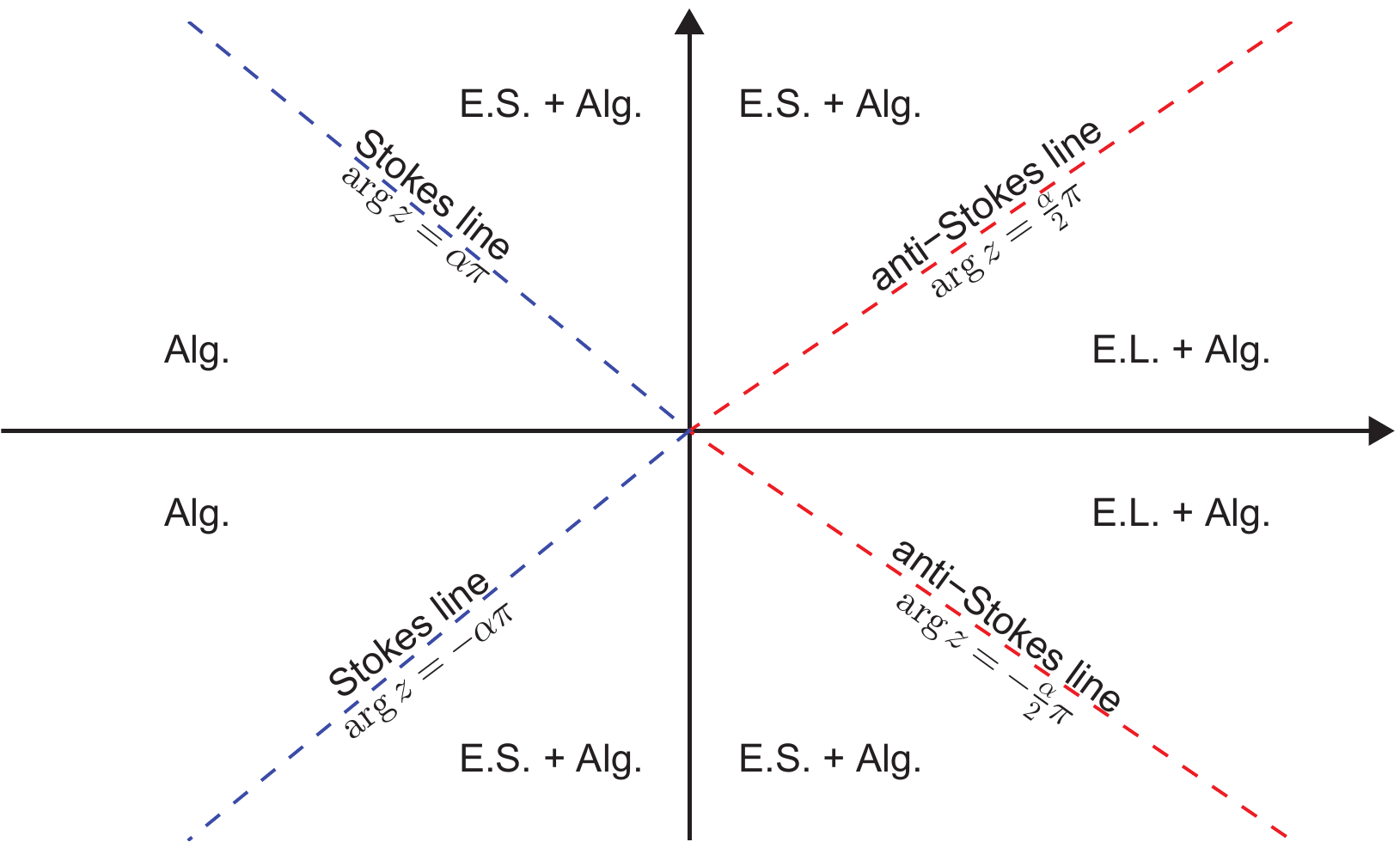} 
\caption{Asymptotic behavior of $E_{\alpha,\beta}^{\gamma}(z)$ in different regions of the complex plane for $0<\alpha<1$. {E.S.}: \emph{exponentially small}, {E.L.}: \emph{exponentially large} and {Alg.}: \emph{algebraic}}
\label{fig:Fig_StokesLine}
\end{figure}

It turns out that the behavior of $E_{\alpha,\beta}^{\gamma}(z)$ for large values of $|z|$ varies from exponential to algebraical depending on the sector of the complex plane where $z$ lies. Furthermore, it is found that the parameter $\alpha$ is the one that (most prominently) controls the asymptotic properties of the Prabhakar function. Specifically, the lines $\arg z=\pm\alpha\pi/2$ are {\em anti-Stokes lines}, where the function changes its behavior from increasing to decreasing. Whereas, $\arg z=\pm \alpha\pi$ are the {\em Stokes lines}, where the exponential term quickly decays leaving just a predominant algebraic term \cite{Paris2019}. This behavior is depicted in Figure \ref{fig:Fig_StokesLine} for $0 < \alpha < 1$ (the acronyms ``{E.S.}'', ``{E.L.}'' and ``{Alg.}'' stand for ``\emph{exponentially small}'', ``\emph{exponentially large}'', and ``\emph{algebraic}'', respectively). When $1<\alpha<2$ the Stokes lines $\arg z = \pm \alpha \pi$ collapse onto the negative real axis and there are no more regions in which the Prabhakar function presents a significant algebraic expansion.

The explicit form of the exponential and algebraic terms, together with an algorithm for the evaluation of their coefficients, were described in \cite{Paris2010,Paris2019}. In particular the algebraic expansion is
\vskip -12pt
\[
	{\mathcal A}_{\alpha,\beta}^{\gamma}(z) = \frac{z^{-\gamma}}{\Gamma(\gamma)} \sum_{k=0}^{\infty} \frac{(-1)^{k} \Gamma(k+\gamma)}{k! \Gamma(\beta-\alpha(k+\gamma))} z^{-k}
\]
\vskip -3pt \noindent
while the exponential expansion reads
\vskip -10pt
\[
	{\mathcal E}_{\alpha,\beta}^{\gamma}(z) = \frac{1}{\Gamma(\gamma)} \eu^{z^{1/\alpha}} z^{\frac{\gamma-\beta}{\alpha}} \frac{1}{\alpha^{\gamma}} \sum_{k=0}^{\infty} c_{k} z^{-\frac{k}{\alpha}} ,
\]
where $c_{k}$ are obtained from the inverse factorial expansion of
\begin{equation}\label{eq:InverseFactorialExpansion}
	F_{\alpha,\beta}^{\gamma}(s) := \frac{\Gamma(\gamma+s)\Gamma(\alpha s + \psi)}{\Gamma(s+1)\Gamma(\alpha s + \beta)}
	= \alpha^{1-\gamma} \left( 1 + \sum_{j=1}^{\infty} \frac{c_{j}}{(\alpha s + \psi )_j} \right)
\end{equation}
for $|s|\to \infty$ in $|\arg(s)| \le \pi - \epsilon$ and any arbitrarily small $\epsilon >0$. As usual, $(x)_j$ denotes the Pochhammer symbol and $\psi=1-\gamma+\beta$.

One can now summarize the asymptotic behavior of the Prabhakar function, for $|z|\to\infty$, as follows (see \cite{Paris2019})
\[
	\begin{aligned}
	0 < \alpha \le 1 \, &: \\
	E_{\alpha,\beta}^{\gamma}(z) &\sim \left\{\begin{array}{ll}
		{\mathcal E}_{\alpha,\beta}^{\gamma}(z) + {\mathcal A}_{\alpha,\beta}^{\gamma}(z \eu^{\mp\pi\iu})
		& |\arg z | < \frac{1}{2}\alpha \pi \\
		{\mathcal A}_{\alpha,\beta}^{\gamma}(z \eu^{\mp\pi\iu}) + {\mathcal E}_{\alpha,\beta}^{\gamma}(z)
		& \frac{1}{2}\alpha \pi < |\arg z | < \alpha \pi \\
		{\mathcal A}_{\alpha,\beta}^{\gamma}(z \eu^{\mp\pi\iu}) & \alpha \pi < |\arg z | \le \pi \\
	\end{array}\right. \\
	1 < \alpha < 2 \, &: \\
	E_{\alpha,\beta}^{\gamma}(z) &\sim \left\{\begin{array}{ll}
		{\mathcal E}_{\alpha,\beta}^{\gamma}(z) + {\mathcal A}_{\alpha,\beta}^{\gamma}(z \eu^{\mp\pi\iu})
		& |\arg z | < \frac{1}{2}\alpha \pi \\
		{\mathcal A}_{\alpha,\beta}^{\gamma}(z \eu^{\mp\pi\iu}) + {\mathcal E}_{\alpha,\beta}^{\gamma}(z) + {\mathcal E}_{\alpha,\beta}^{\gamma}(z \eu^{\mp 2\pi\iu})
		& \frac{1}{2}\alpha \pi < |\arg z | < \pi \\
	\end{array}\right.
	\end{aligned}
\]
where the sign in $\eu^{\mp\pi\iu}$ and $\eu^{\mp 2\pi\iu}$ is negative when $z$ lies in the upper complex half-plane and positive otherwise. Note that some of the cases discussed above differ just by the order in which ${\mathcal E}_{\alpha,\beta}^{\gamma}(z)$ and ${\mathcal A}_{\alpha,\beta}^{\gamma}(z \eu^{\mp\pi\iu})$ appear. The reason for this choice comes from the fact that, for the sake of clarity, we prefer to give priority in the expression to the dominant term. The term ${\mathcal E}_{\alpha,\beta}^{\gamma}(z \eu^{\mp 2\pi\iu})$ can be usually neglected, except when $\arg z$ is close to $\pi$ since in this case it becomes comparable to ${\mathcal E}_{\alpha,\beta}^{\gamma}(z)$.

	To complete the analysis we consider the case $\alpha \ge 2$. In detail, one finds
\vskip -13pt
\[
	\alpha \ge 2 \, : \quad
	E_{\alpha,\beta}^{\gamma}(z) \sim
	\sum_{r=-P}^{P} {\mathcal E}_{\alpha,\beta}^{\gamma}(z\eu^{2\pi\iu r}) + {\mathcal A}_{\alpha,\beta}^{\gamma}(z \eu^{\mp\pi\iu}) \, , \,\,\,  \quad |\arg z | \le \pi \, , \\
\]
\vskip -4pt \noindent
where $P$ is an integer number such that $2P+1$ is the smallest odd integer satisfying $2P+1>1/\alpha$. For the latter
it is found that the anti-Stokes lines collapse onto the negative real axis, thus implying that the asymptotic behavior is always exponentially large. Indeed, we included the algebraic term just for the sake of completeness, though it turns out to be always negligible.

Providing an explicit formulation of the coefficients $c_{k}$ in the exponential expansion ${\mathcal E}_{\alpha,\beta}^{\gamma}(z)$ is not an easy task and a sophisticated algorithm is discussed in \cite{Paris2010}. The first few values of the $c_k$'s have however been explicitly computed in \cite{GarraGarrappa2018} and \cite{Paris2019}:
\[
	\begin{aligned}
	c_{0} &= 1 \\
	c_{1} &= \frac{(\gamma-1)}{2} \left( \alpha \gamma + \gamma - 2 \beta\right) \\
	c_{2} &= \frac{(\gamma-1)(\gamma-2)}{24} \left[ 3(\alpha+1)^2\gamma^2 - (\alpha+1)(\alpha+12\beta+5) \gamma + 12\beta(1+\beta) \right] \\
	c_{3} &= \frac{(\gamma-1)(\gamma-2)(\gamma-3)}{48} \bigl[ \gamma^3(1+\alpha)^3 - \gamma^2(1+\alpha)^2(5+\alpha+6\beta)\\
	      & \phantom{\qquad} +2\gamma(1+\alpha)(3+\alpha(1+\beta)+11\beta+6\beta^2) - 8\beta(1+\beta)(2+\beta)\bigr] \, . \\
	\end{aligned}
\]
Showing a longer list of these coefficients is not particularly instructive, hence we just refer again to the algorithm provided in \cite{Paris2010}.

Specific results for the asymptotic expansion along the negative real semi-axis, which is of most interest for the study of relaxation phenomena, have been provided, for instance, in \cite{MainardiGarrappa2015_JCP}. In this case, focusing just on the leading term of the expansion, one can infer that
	\begin{equation}
	E^{\gamma}_{\alpha, \beta} (-t^\alpha) \sim
	\begin{cases}
	\frac{1}{\Gamma(\beta-\alpha\gamma)}t^{-\alpha \gamma} \quad & \beta \neq \alpha\gamma\\
    -\frac{\gamma}{\Gamma(-\alpha)}t^{-\alpha \gamma-\alpha} \quad & \beta = \alpha\gamma
	\end{cases}
\end{equation}
as $t\to+\infty$. 	

\vspace*{-2pt}

\section{Prabhakar fractional calculus} \label{sec:pracalculus}

\setcounter{section}{5} \setcounter{equation}{0} 
	
After an extensive review of the main properties of the Prabhakar function, we are now ready to introduce and discuss a generalization of fractional calculus based on this special function, that we shall call {\em Prabhakar fractional calculus}.

One of the main contribution of T. R. Prabhakar \cite{Prabhakar1971} consisted in the formulation of a new class of weakly-singular linear Volterra operators based on the function
\vskip -10pt
\begin{align} \label{eq:PrabhakarKernel}
	e^{\gamma}_{\alpha,\beta}(t; \lambda) = t^{\beta-1}E^{\gamma}_{\alpha,\beta}\left( \lambda t^{\alpha}\right),
	\quad \alpha >0 \, ,
\end{align}
\vskip -3pt \noindent
known in the literature as the {\em Prabhakar kernel}, with $\lambda$ being a real (or complex) number. Most notably, one finds that for $\gamma = 0$ or $\lambda=0$ this kernel reduces to
\vskip -12pt
\begin{equation}\label{eq:SpecialCasesPrabhakar}
	e^{0}_{\alpha,\beta}(t;\lambda) = e^{\gamma}_{\alpha,\beta}(t;0) = \frac{t^{\beta-1}}{\Gamma(\beta)} \, .
\end{equation}

	Let $f \in L^1[t_0,T]$, then the integral operator originally proposed in \cite{Prabhakar1971}, which is nowadays universally recognized as the {\em Prabhakar fractional integral}, reads
\vskip -10pt
\begin{equation}\label{eq:PrabhakarIntegral}				
\begin{aligned}
	\JPr^{\gamma}_{\alpha,\beta, \lambda; t_0}f(t)
	&\coloneqq \bigl( e^{\gamma}_{\alpha,\beta}( \, \cdot \, ; \lambda) \ast f \bigr) \\
	&=\int_{t_0}^{t}(t-u)^{\beta-1}E^{\gamma}_{\alpha,\beta}\left[\lambda (t-u)^{\alpha} %
	\right]f(u) \du u,  \quad \alpha,\beta>0,
\end{aligned}
\end{equation}
where the condition $\alpha>0$ is necessary for the convergence of the series expansion (\ref{eq:PrabhakarFunction}), whereas the restriction $\beta>0$ is imposed in order to guarantee the convergence of the integral. In view of \eqref{eq:SpecialCasesPrabhakar}, we have that the Prabhakar kernel in \eqref{eq:PrabhakarIntegral} reduces to the Gel'fand--Shilov distribution \cite{Giusti2019}, which implies that
\vskip -10pt
\be
\JPr^{0}_{\alpha,\beta,\lambda;t_0} f(t) = \JPr^{\gamma}_{\alpha,\beta,0;t_0} f(t) = J^{\beta}_{t_0} f(t) \, ,
\ee
\ie \eqref{eq:PrabhakarIntegral} reduces to the RL integral for $\gamma = 0$ or $\lambda=0$. Besides, another interesting relation between the RL and Prabhakar integrals is related to the {\em series representation} of the latter. Indeed, let $f \in L^1[t_0,T]$ and $\alpha, \beta > 0$, then it is easy to see that (see \cite{Giusti2018_NODY})
\vskip -10pt
\be \label{eq:praseriesint}
\JPr^{\gamma}_{\alpha,\beta, \lambda; t_0}f(t) =
\sum _{k=0} ^\infty \frac{(\gamma)_k \lambda^k}{k!} \, J^{\alpha k + \beta}_{t_0} f(t) \, .
\ee

Additionally, a property that finds several applications when dealing with Prabhakar's theory is given by
\vskip -10pt
\be \label{eq:PrabIntSum}	
\qquad
\JPr^{\gamma}_{\alpha, \beta, \lambda; t_0}
	\Bigl[(t-t_0)^{\mu-1}E_{\alpha ,\mu}^{\sigma}\big( \lambda (t-t_0)^{\alpha}\big)\Bigr]
	= (t-t_0)^{\beta+\mu-1} E_{\alpha ,\beta+\mu}^{\gamma+\sigma} \big( \lambda (t-t_0)^{\alpha}\big)
\ee
or, in other terms,
\vskip -10pt
\begin{equation}
	\JPr^{\gamma}_{\alpha, \beta, \lambda; t_0}
	\Bigl[e^{\sigma}_{\alpha,\mu}( t - t_0 \, ; \lambda) \Bigr]
	= e_{\alpha ,\beta+\mu}^{\gamma+\sigma} ( t - t_0 \, ; \lambda) \, .
\end{equation}

	The Prabhakar fractional integral has been extensively studied in \cite{KilbasSaigoSaxena2004}, where it was shown that
$\JPr^{\gamma}_{\alpha,\beta, \lambda; t_0}: L_1[t_0,T] \to L_1[t_0,T]$ is a bounded operator (see \cite[Theorem 4]{KilbasSaigoSaxena2004}). Most importantly, in \cite{KilbasSaigoSaxena2004} a family of left-inverse operators of the Prabhakar fractional integral was identified, \ie
\begin{equation}\label{eq:PrabhakarDerivativeRLGeneral}
	\bm{\mathcal{D}}^{\gamma}_{\alpha, \beta, \lambda; t_0} f(t)
	       = \DR^{\mu}_{t_0} \JPr^{-\gamma}_{\alpha, \mu-\beta,\lambda; t_0}f(t)
	\quad \forall \mu \in \R \, \, \text{with} \, \, \mu >\beta
\end{equation}
with $f \in L_1[t_0,T]$. The presence of the fractional-order RL derivative in (\ref{eq:PrabhakarDerivativeRLGeneral}) makes this family of operators too convoluted for most practical applications. Hence, in most of the subsequent works (\emph{e.g.}, see \cite{GarraGorenfloPolitoTomovski2014,GarrappaMainardiMaione2016,PolitoTomovski2016}) the special case $\mu = m = \left\lceil \beta \right\rceil$ is usually preferred. This then leads to what is now commonly known as the {\em Prabhakar derivative} of RL type, \ie
\begin{equation}\label{eq:PrabhakarDerivativeRL}
\begin{aligned}
	\DPrRL^{\gamma}_{\alpha, \beta, \lambda; t_0}f(t)
	&\coloneqq D^{m} \, \JPr^{-\gamma}_{\alpha, m-\beta, \lambda; t_0} f(t)  \\
	&= \frac{\du^m}{\du t^m}
	\int_{t_0}^{t}(t-u)^{m-\beta-1}E^{-\gamma}_{\alpha,m-\beta}\left[\lambda (t-u)^{\alpha} 	\right]f(u) \du u
	\, .
\end{aligned}
\end{equation}

	In analogy with the classical theory of fractional differential operators, one can introduce a regularized version of the Prabhakar derivative (\ref{eq:PrabhakarDerivativeRL}). This is achieved by exchanging the order according to which the operations of integration and differentiation appear in \eqref{eq:PrabhakarDerivativeRL}. This procedure is known as the {\em Caputo-like regularization} of the fractional Prabhakar derivative and it was first performed by M. D'Ovidio and F. Polito in \cite{DOvidioPolito2018}. Then, let $f \in AC^m [t_0,T]$, the {\em regularized Prabhakar derivative} is defined as
\begin{equation}\label{eq:PrabhakarDerivativeC}
\begin{aligned}
	\DPrC^{\gamma}_{\alpha, \beta, \lambda; t_0}f(t)
	&\coloneqq \JPr^{-\gamma}_{\alpha, m-\beta, \lambda; t_0} D^{m}  f(t)  \\
	&= 	\int_{t_0}^{t}(t-u)^{m-\beta-1}E^{-\gamma}_{\alpha,m-\beta}\left( \lambda (t-u)^{\alpha} \right)  f^{(m)}(u) \, 		\du u .
\end{aligned}
\end{equation}

This operator turns out to act as another left-inverse of the Prabhakar fractional integral, as shown in \cite{Giusti2019}.

One reason why \eqref{eq:PrabhakarDerivativeC} is understood as a regularization of \eqref{eq:PrabhakarDerivativeRL} is due to the fact that this procedure makes the solutions to the eigenvalue problem associated to \eqref{eq:PrabhakarDerivativeC} non-singular at $t=t_0$. Property which is not shared by the solutions to the eigenvalue problem associated to the (RL-type) Prabhakar derivative. Furthermore, again differently from the RL-type case, one can set up a well-posed Cauchy problem for the regularized Prabhakar derivative with initial conditions involving solely ordinary derivatives of the unknown function at $t=t_0$. Besides, one also has that
\[
	\DPrC^{\gamma}_{\alpha, \beta, \lambda; t_0} t^{k} = 0
	, \quad k = 0,1,\dots,m-1, \quad m = \left\lceil \beta \right\rceil \, .
\]

	It is now easy to see that for $\gamma = 0$ or $\lambda=0$ (\ref{eq:PrabhakarDerivativeRL}) and \eqref{eq:PrabhakarDerivativeC} reduce to the standard RL and Caputo derivatives, namely \eqref{eq:RL_Derivative} and \eqref{eq:C_Derivative}. Moreover, it turns out that \eqref{eq:PrabhakarDerivativeRL} and \eqref{eq:PrabhakarDerivativeC} satisfy
\begin{align} \label{eq:equivalence}
	\DPrC^{\gamma}_{\alpha, \beta, \lambda; t_0} f(t) =
	\DPrRL^{\gamma}_{\alpha, \beta, \lambda; t_0}
	\left[ f(t) - \sum_{k=0}^{m-1}\frac{(t-t_0)^k}{k!}f^{(k)}(t_{0}^+) \right] \, ,
\end{align}
further strengthening their similarity to canonical fractional derivatives. Furthermore, using this last property A. Giusti \cite{Giusti2019} was then able to frame (to some extent) Prabhakar's theory within the scheme of Kochubei's {\em general fractional calculus} \cite{Kochubei2011}. Note that, to this regard, some preliminary analysis had been carried out in \cite{GorskaKochubei} for a very specific problem.

	Before concluding this general summary of the main definitions and properties of the Prabhakar fractional integral and derivatives, it is worth stressing that a similar result to the series expansion in \eqref{eq:praseriesint} can be trivially computed also for the operators \eqref{eq:PrabhakarDerivativeRL} and \eqref{eq:PrabhakarDerivativeC}.

\subsection{Laplace domain} 

It is now of paramount importance to recall a few results concerning the effects of integral transforms on Prabhakar operators.

First, let us consider the Laplace transform of the Prabhakar fractional integral, with $t_0 = 0$. Since \eqref{eq:PrabhakarIntegral} is a convolution-type integral, in the Laplace sense, it is easy to see that
\vskip -12pt
\be
\mathcal{L} \Big[ \JPr^{\gamma}_{\alpha, \beta, \lambda; 0} f(t) \, ; \, s \Big]
=
\mathcal{L} \Big[ e^{\gamma}_{\alpha,\beta}( t \, ; \lambda) \, ; \, s \Big] \, \wt{f} (s)
=
\frac{s^{\alpha\gamma - \beta}}{\bigl(s^{\alpha} - \lambda \bigr)^{\gamma}}
\, \wt{f} (s) \, ,
\ee
\vskip -3pt \noindent
where in the last step we used the result in \eqref{ML3_LT}, with $\real(s) > 0$ and $|s|>|\lambda|^{\frac{1}{\alpha}}$. Similarly, for the Prabhakar derivative \eqref{eq:PrabhakarDerivativeRL} one finds
\vskip -10pt
\be
\notag \mathcal{L} \Big[ \DPrRL^{\gamma}_{\alpha, \beta, \lambda; 0} f(t) \, ; \, s \Big]
&=&
s^{\beta-\alpha\gamma} (s^\alpha - \lambda)^\gamma \, \wt{f} (s)\\
& &-
\sum _{k=0} ^{m-1}
s^{m-k-1} \, \Big( \JPr^{-\gamma}_{\alpha, m - \beta - k, \lambda; 0} f  \Big) (0^+) \, ,
\ee
\vskip -3pt \noindent
with $m = \ceil{\beta}$, whose proof is simply based on the LT of the integer-order derivative, the result in \eqref{ML3_LT}, and that
(see \cite{KilbasSaigoSaxena2004})
\vskip -10pt
$$
\DR _{t_0} ^\sigma \JPr^{\gamma}_{\alpha, \beta, \lambda; t_0} f(t) =
\JPr^{\gamma}_{\alpha, \beta - \sigma, \lambda; t_0} f(t) \, ,
$$
for $f \in L_1 [t_0, T]$, $\sigma > 0$ and $\beta > \sigma$, and in particular $D^k \JPr^{\gamma}_{\alpha, \beta, \lambda; t_0} f(t) =
\JPr^{\gamma}_{\alpha, \beta - k, \lambda; t_0} f(t)$ with $k$ being an integer. For the sake of completeness, one also has that (see \cite{KilbasSaigoSaxena2004})
\vskip -11pt
$$
J_{t_0} ^\sigma \JPr^{\gamma}_{\alpha, \beta, \lambda; t_0} f(t) =
\JPr^{\gamma}_{\alpha, \beta + \sigma, \lambda; t_0} f(t) =
\JPr^{\gamma}_{\alpha, \beta, \lambda; t_0} \Big[
J_{t_0} ^\sigma f(t) \Big] \, ,
$$
\vskip -3pt \noindent
for $f \in L_1 [t_0, T]$ and $\sigma > 0$.

	Analogously, one also finds
\vskip -12pt
	\be
	\notag \mathcal{L} \Big[ \DPrC^{\gamma}_{\alpha, \beta, \lambda; 0} f(t) \, ; \, s \Big] 	
	=
	s^{\beta-\alpha\gamma} (s^\alpha - \lambda)^\gamma \Bigg[ \wt{f} (s)
	-
	\sum _{k=0} ^{m-1}
	s^{-k-1} \, f^{(k)} (0^+)
	\Bigg]
	\, ,
	\ee
\vskip -4pt \noindent
where $m = \ceil{\beta}$, as usual.

\subsection{Eigenvalue problem} 

Let us consider the Cauchy problem
\be \label{IVP-Caputo}
\left\{
\begin{aligned}
& \DPrC^{\gamma}_{\alpha, \beta, \lambda; 0} y(t) = A \, y(t) \, , \\
& y (0^+) = \xi_0 \, , \,\, y' (0^+) = \xi_1 \, , \ldots , \,\, y^{(m-1)} (0^+) = \xi_{m-1} \, ,
\end{aligned}
\right.
\ee
with $\alpha, \beta >0$, $m = \ceil{\beta}$, and $A, \xi_0, \ldots, \xi_{m-1} \in \R$.

The Cauchy problem in \eqref{IVP-Caputo} is well-posed and one can compute its solution by means of the Laplace transform method. Specifically, one finds
\begin{equation}\label{eq:SolutionPrabTestEquation}
y(t) = \sum_{j=0}^{m-1} \sum_{k=0}^{\infty} A^k t^{\beta k + j } E_{\alpha,\beta k + j + 1}^{\gamma k}(\lambda t^{\alpha}) \xi_j
\end{equation}
which thus gives the eigenfunction of $\DPrC^{\gamma}_{\alpha, \beta, \lambda; 0}$ corresponding to the eigenvalue $A$. Note that, if $\gamma=0$ or $\lambda=0$, then
\be
y(t) = \sum_{j=0}^{m-1} \sum_{k=0}^{\infty}  \frac{A^k t^{\beta k + j }}{\Gamma(\beta k + j +1 )} \, \xi_j =
\sum_{j=0}^{m-1} t^{j} E_{\beta,j+1}(A t^{\beta} ) \, \xi_j \, ,
\ee
which is indeed the eigenfuction of the Caputo derivative $\DC^{\beta}_0$ corresponding to the eigenvalue $A$ (see \cite{GorenfloKilbasMainardiRogosin2014}), as expected.

Moving on to the eigenvalue problem for the (RL-type) Prabhakar derivative, let us consider the Cauchy problem
\begin{equation}
\left\{ \begin{array}{l}
        \DPrRL^{\gamma}_{\alpha, \beta, \lambda;0} y(t) = A y(t) \\
        \Bigl(\JPr^{-\gamma}_{\alpha, m-\beta-j, \lambda; 0} y \Bigr) (0^+) =
\xi_j , \quad j = 0, 1, \dots, m-1 \, . \\
\end{array} \right.
\end{equation}
Following a procedure akin to the one presented above we find that the solution to this problem reads
\vskip -10pt
\begin{equation}\label{eq:SolutionPrabRLTestEquation}
y(t) = \sum_{j=0}^{m-1} \sum_{k=0}^{\infty}  A^k t^{\beta k +\beta - j
- 1} E_{\alpha,\beta k + \beta - j }^{\gamma (k+1)}(\lambda \, t^{\alpha}) \,
\xi_j
\end{equation}
and hence the functions
\[
\sum_{k=0}^{\infty}  A^k t^{\beta k +\beta - j - 1} E_{\alpha,\beta k
+ \beta - j }^{\gamma (k+1)}( \lambda \, t^{\alpha})
  , \quad j=0,1,\dots,\left\lceil \beta\right\rceil-1
\]
are eigenfunctions of $\DPrRL^{\gamma}_{\alpha, \beta, \lambda; 0}$. Besides, for
$\gamma=0$ or $\lambda=0$, in light of
(\ref{eq:SpecialCasesPrabhakar}) we have that Eq.
(\ref{eq:SolutionPrabRLTestEquation}) becomes
\[
        y(t) = \sum_{j=0}^{m-1} \sum_{k=0}^{\infty}  \frac{ A^k t^{\beta k +
\beta - j -1 }}{\Gamma(\beta k + \beta - j )} \, \xi_j =
\sum_{j=0}^{m-1} t^{\beta-j-1} E_{\beta,\beta-j}(A \, t^{\beta}) \, \xi_j
\]
thus yielding the known eigenfunctions $t^{\beta-j-1}
E_{\beta,\beta-j}(A \, t^{\beta})$, $j=0,1,\dots,\left\lceil
\beta\right\rceil-1$ of the RL fractional derivative $\DR^{\beta}_0$
\cite{GorenfloKilbasMainardiRogosin2014}. In particular, note that for $j=m-1$ one has that $t^{\beta-j-1} = t^{\beta-\ceil{\beta}}$ which carries a weakly singular behavior if $\beta \notin \N$ since in this case $\beta-\ceil{\beta} < 0$.

\subsection{Operators of Gr\"{u}nwald-Letnikov type} 

It is now worth mentioning that a Gr\"{u}nwald--Letnikov (GL) type formulation of the Prabhakar operators was proposed by R. Garrappa in \cite{GarrappaHN2016}. Originally, this interpretation of the theory was conceived for the particular case $\beta=\alpha\gamma$ that turns out to be involved in the time-domain representation of the Havriliak--Negami relaxation. Its generalization to $\beta>0$ is however immediate and is achieved by following the same steps illustrated in \cite{GarrappaHN2016} and by taking profit of Lubich's theory \cite{Lubich1988a,Lubich1988b,Lubich2004} of generalized convolution quadratures. Indeed, one can find that, under suitable conditions, the Prabhakar integral \eqref{eq:PrabhakarIntegral} is equivalent to the GL integral
\vskip -12pt
\begin{equation}\label{eq:GL_Prab_Int}
	\bm\bar{\JPr}^{\gamma}_{\alpha, \beta, \lambda; t_0} f(t) = \lim_{h\to 0} \frac{h^{\beta}}{(1-h^{\alpha}\lambda)^{\gamma}} \sum_{j=0}^{\left\lfloor \frac{t-t_0}{h}\right\rfloor} W^{(-)}_{n-j} f(t-jh) ,
\end{equation}
\vskip -4pt \noindent
where the coefficients $W^{(-)}_{n}$ depend on $\alpha$, $\beta$, $\gamma$, $\lambda$ and $h$. For $\gamma\not=0$, $W^{(-)}_{n}$ can be evaluated recursively as
\vskip -12pt
\[
	W^{(-)}_0 = 1 , \quad
	W^{(-)}_k = \sum_{j=1}^{k} \left( \frac{(1-\gamma)j}{k}-1 \right) \bar{\omega}_{j} W^{(-)}_{k-j}	,
\]
\vskip -4pt \noindent
with
\vskip -11pt
\[
\bar{\omega}_{j} = \frac{\omega_{j}^{(\beta/\gamma)} - h^{\alpha}\lambda\omega_{j}^{(\beta/\gamma-\alpha)}}{1-h^{\alpha} \lambda} \, ,
\]
while the $\omega_{j}^{(\alpha)}$ are related to the classical binomial coefficients according to
\vskip -10pt
\[
\omega_{j}^{(\alpha)} = (-1)^j \binom{\alpha}{j} = (-1)^j \frac{\Gamma(\alpha+1)}{\Gamma(j+1) \Gamma(\alpha-j+1)}\, .
\]	

Clearly, the GL--Prabhakar integral in \eqref{eq:GL_Prab_Int} generalizes the classical GL operator corresponding to the RL integral \eqref{eq:RL_Integral}. The RL integral is instead recovered from \eqref{eq:GL_Prab_Int} by setting $\lambda=0$ since, in this case, one has $W^{(-)}_j=\omega_{j}^{(-\beta)}$.

In a similar fashion, a GL operator corresponding to Prabhakar derivative $\DPrRL^{\gamma}_{\alpha, \beta, \lambda; t_0}$ can also be obtained as
\vskip -12pt
\[
	{}^{\texttiny{RL}}\!\bm\bar{\bm{\mathcal{D}}}^{\gamma}_{\alpha, \beta, \lambda; t_0} f(t) = \lim_{h\to 0} \frac{(1-h^{\alpha}\lambda)^{\gamma}}{h^{\beta}} \sum_{j=0}^{\left\lfloor \frac{t-t_0}{h}\right\rfloor} W^{(+)}_{n-j} f(t-jh) ,
\]
\vskip -4pt \noindent
with
\vskip -12pt
\[
	W^{(+)}_0 = 1 , \quad
	W^{(+)}_k = \sum_{j=1}^{k} \left( \frac{(1+\gamma)j}{k}-1 \right) \bar{\omega}_{j} W^{(+)}_{k-j}	,
\]
where the coefficients $\bar{\omega}_{j}$ are the same as the one introduced above. It is worth noting, again, that the classical GL derivative, equivalent to the RL derivative \eqref{eq:RL_Derivative}, is recovered for $\lambda=0$.

	Finally, in light of \eqref{eq:equivalence} a GL operator related to the regularized Prabhakar derivative $\DPrC^{\gamma}_{\alpha, \beta, \lambda; t_0}$ can be easily obtained. The importance of these operators lies in their peculiar proneness toward their use as numerical approximation of the Prabhakar operators, once the step-size $h>0$ has been fixed.

\vspace*{-4pt}
\section{Physical applications} \label{sec:physics}

\setcounter{section}{6} \setcounter{equation}{0} 

	The aim of this section is to briefly summarize the main results and physical implications of Prabhakar calculus. Of course, the literature on this subject is rather vast and, at time, dispersive, hence our plan is to offer a general overview of the key physical features of these models that fall in the realm of \emph{anomalous phenomena}. Specifically, we will start off by describing how this new framework can describe dielectric relaxation processes such as the one in the Davidson--Cole and the Havriliak--Negami models. Then, we will show how the electro-mechanical correspondence theorized by B. Gross and R. M. Fuoss paves the way to a brand new class of viscoelastic models. Finally, we summarize the main results concerning the role of the Prabhakar function (and the corresponding operators) in the theory of anomalous diffusion on complex media.

\subsection{Anomalous relaxation in dielectrics} \label{subsec:dielectrics} 

It is fairly well known in the literature \cite{CC1, CC2, DC, Vaughan1979, Debye, WangCD2019, WangHN2019} that materials that display glass-liquid transitions or amorphous polymers tend to deviate substantially from the standard Debye relaxation model.	 Setting aside for the moment the Cole--Cole relaxation, that can be described in terms of standard fractional calculus \cite{CapelasMainardiVaz2011, GarrappaMainardiMaione2016, GorskaPLA2019}, let us focus on the Havriliak--Negami model. First, recalling that we are dealing with a causal linear system, we can recast Eq. \eqref{HN} in the Laplace domain, thus the response function will read \cite{CapelasMainardiVaz2011, GarrappaMainardiMaione2016, GorskaJPA2018}
\vskip -12pt
	\be \label{HN-Laplace}
	\wt{\phi} _{\rm HN} (s) =
	\wt{\chi} _{\rm HN} (s) = \frac{1}{[1 + (s \, \tau _{\star})^\alpha]^\gamma} 		\, ,
	\quad 0 < \alpha \leq 1 \, , \,\,\, \gamma > 0 \, ,
	\ee
	then, taking profit of Eq. \eqref{ML3_LT}, one gets \cite{GarrappaMainardiMaione2016}
	\be
	\phi_{\rm HN} (t) =
	\frac{1}{\tau _{\star}} \left( \frac{t}{\tau _{\star}} \right) ^{\alpha \gamma - 1}
	E ^\gamma _{\alpha, \alpha \gamma} \left[ - \left( \frac{t}{\tau _{\star}} \right) ^{\alpha}
	\right] \, ,
	\ee
	that, in view of Eq. \eqref{relaxation-function}, leads to the relaxation function \cite{GarrappaMainardiMaione2016}
	\be \label{HN-relaxation}
	\Psi_{\rm HN} (t) =
	1 - \left( \frac{t}{\tau _{\star}} \right) ^{\alpha \gamma }
	E ^\gamma _{\alpha, \alpha \gamma+1} \left[ - \left( \frac{t}{\tau _{\star}} \right) ^{\alpha}
	\right] \, .
	\ee	

	A few comments are now in order. First of all, it is worth mentioning the asymptotic behaviors of $\Psi_{\rm HN} (t)$, namely
	\be
	\Psi_{\rm HN} (t) \sim
	\left\{
	\begin{aligned}
	1 - \frac{(t / \tau _\star)^{\alpha \gamma}}{\Gamma (1 + \alpha \gamma)} \, , \qquad 	 	
	&\mbox{for} \,\, t \ll \tau _\star \, , \\
	\frac{\gamma (t / \tau _\star)^{-\alpha}}{\Gamma (1 - \alpha)} \, , \qquad
	&\mbox{for} \,\, t \gg 	\tau _\star \, .
	\end{aligned}
	\right.
	\ee
	Besides, it is also straightforward to prove that the \emph{relaxation equation} for the Havriliak-Negami model is given by
	\be
	\DPrC ^{\gamma} _{\alpha , \alpha \gamma , \tau _\star^{-\alpha}; \, 0} \Psi_{\rm HN} (t) =
	- \frac{1}{\tau _\star^{\alpha \gamma}} \, , \qquad \Psi_{\rm HN} (0) = 1 \, ,
	\ee
	highlighting the strict connection between this model and Prabhakar calculus.
	
	Finally, one can easily infer that the Davidson--Cole and the Cole--Cole models are particular realizations of the Havriliak--Negami relaxation. Indeed, the Davidson--Cole dielectric response is obtained from the Havriliak--Negami one by setting $\alpha = 1$ in Eq. \eqref{HN}, leading to a relaxation function that reads \cite{GarrappaMainardiMaione2016}
\vskip -12pt
	\be
	\Psi _{\rm DC} (t) = 1 - \left( \frac{t}{\tau _\star} \right) ^\gamma
	E^\gamma _{1, \gamma +1} \left(- \frac{t}{\tau _\star} \right) =
	\frac{\Gamma (\gamma, t/\tau _\star)}{\Gamma (\gamma)} \, ,
	\ee
	with $\Gamma (a, z) = \int_z ^\infty u^{a - 1} \, e^{-u} \, \d u$ denoting the incomplete gamma function, and a relaxation equation
\vskip -10pt
	\be
	\DPrC ^{\gamma} _{1, \gamma , \tau _\star^{-1}; \, 0} \Psi_{\rm DC} (t)
	=
	- \frac{1}{\tau _\star^{\gamma}} \, , \qquad \Psi_{\rm DC} (0) = 1 \, .
	\ee
	Analogously, it is easy to see that the Cole--Cole model can be derived from Eq. \eqref{HN} by setting $\gamma = 1$. This, in turn, leads the relaxation function in Eq. \eqref{psi-CC} that corresponds to an evolution equation governed by Caputo derivatives \cite{GarrappaMainardiMaione2016, GorskaPLA2019}, \ie
\vskip -10pt
	\be
	\DC_{0} ^\alpha \Psi _{\rm CC} (t) = - \frac{1}{\tau _{\star} ^\alpha} \, \Psi _{\rm CC} (t)
	\, ,
	\qquad \Psi_{\rm CC} (0) = 1 \, ,
	\ee
	which is indeed the canonical fractionalization of the standard relaxation problem, namely
	$\Psi ' (t) = - (1/\tau) \, \Psi (t)$ with $\Psi (0) = 1$.
	
	Before concluding this Section we wish to stress that the models discussed here clearly do not exhaust the full realm of theoretical set-ups that have been proposed, over the years, to describe the features of anomalous dielectrics. In this work we simply report the most known modifications of Debye's theory that display a strict connection with the Prabhakar function and calculus. For a more detailed analysis of the literature on anomalous dielectric relaxation and its relation to the notion of complete monotonicity, we refer the interested reader to \cite{CapelasMainardiVaz2011, MainardiGarrappa2015_JCP, GarrappaMainardiMaione2016}.
	
	\subsection{Linear Viscoelasticity} 

	It is rather well known in the literature (see \eg \cite{CaputoMainardi1971-1, CaputoMainardi1971-2, RogosinMainardi-2014, GiustiMainardi-2016, Bessel-2017, Bessel-2017-2} and references therein) that fractional calculus and the theory of completely monotone functions play a pivotal role and lie at heart of the modern mathematical formulation of linear viscoelasticity. To this regard, we refer the readers to \cite{Mainardi2010, MainardiSpada2011, Pipkin, GurtinSternberg1962, WaveMotion-2017} for a detailed description of both ordinary and fractional viscoelasticity.
	
	Beside the clear effectiveness of Prabhakar calculus in serving as a fundamental tool for modeling anomalous dielectrics, this novel mathematical structure plays an interesting role also in fractional viscoelasticity. As discussed in \cite{GC-CNSNS2018}, if one recalls the formal duality between viscoelastic models and electrical systems, originally introduced by Gross and Fuoss in \cite{Gross-Fuoss-1, Gross-Fuoss-2, Gross-Fuoss-3} and then revisited in \cite{Gross-Fuoss-GM}, it is easy to see that the anomalous relaxation processes highlighted in Section \ref{subsec:dielectrics} have an analogous representation in linear viscoelasticity. It is therefore natural to wonder about what might happen if one writes the constitutive equation for a certain material in terms of Prabhakar derivatives.
	
	Let us begin with a brief recap of the basics of the linear theory of viscoelasticity. Denoting by $\mathcal{H} ^{N}$ the  Heaviside class, namely the set of causal functions $f(t)$ such that $f \in C^N (\R^+)$ with $N \in \N$, then we denote by $\sigma \in \mathcal{H} ^{N}$ and $\varepsilon \in \mathcal{H} ^{N}$ the uniaxial {\em stress} and {\em strain} functions for a given material, respectively. The key physical information about the system are then stored in the stress-strain relation, which is known as the {\em constitutive equation} for the material, and it can be represented equivalently in differential form or as an integral equation. In particular, the distinctive feature of the integral form is that it relates the stress and the strain through either the {\em creep compliance} $J(t)$ or the {\em relaxation modulus} $G(t)$. In detail, if one considers a quiescent viscoelastic body for $t<0$ and assumes some sufficiently well behaved causal histories, then the general {\em integral form of the constitutive equation} reads
	\be
	\varepsilon (t) = J(0^+) \, \sigma (t) + \int _0 ^t J (t-\tau) \, \sigma ' (\tau) \, \d \tau \, ,
	\ee
\vskip -4pt \noindent
	or equivalently,
\vskip- 12pt
	\be
	\sigma (t) = G(0^+) \, \varepsilon (t) + \int _0 ^t G (t-\tau) \, \varepsilon ' (\tau) \, \d \tau \, .
	\ee
	Since the latter are convolution-type equations, they can be easily manipulated via the standard Laplace transform method. Specifically, these equations can be recast, in the Laplace domain, as
	\be
	\wt{\varepsilon} (s) = s \, \wt{J} (s) \, \wt{\sigma} (s) \, , \quad \wt{\sigma} (s) = s \, \wt{G} (s) \, \wt{\varepsilon} (s)
	\, ,
	\ee
	from which one can infer the so called {\em reciprocity relation}
\vskip -10pt
	\be
	 s \, \wt{J} (s) = \frac{1}{s \, \wt{G} (s)} \, .
	\ee
	Furthermore, it is useful to recall that the creep compliance is a non-negative and non-decreasing function of time on $t\geq 0$, whereas the relaxation modulus is a non-negative and non-increasing function of time on the same domain.
	
	 All that being said, in \cite{GC-CNSNS2018} it was proposed a model representing a straightforward generalization of the renowned {\em Maxwell model} of linear viscoelasticity, where the ordinary derivatives were replaced with Prabhakar's one. In detail, this new proposal is customarily referred to as {\em Maxwell--Prabhakar model} (or, alternatively, as the {\em Giusti--Colombaro model}) of linear viscoelasticity and it is described by the constitutive equation
	 \be \label{GC-model}
	 \sigma (t) + a \, \DPrC ^{\gamma} _{\alpha , \beta, \lambda ; \, 0} \sigma (t)
	 =
	 b \, \DPrC ^{\gamma} _{\alpha , \beta, \lambda ; \, 0} \varepsilon (t)
	 \ee
	 where $\alpha > 0$, $0 < \beta < 1$ and $\gamma , \lambda \in \R$, in general. Besides, since $[\sigma] = \mbox{Force} \cdot \mbox{Area}^{-1}$ and $[\varepsilon] = 1$, one has that $[a] = \mbox{time}^{\beta}$, $[b] = \mbox{Force} \cdot \mbox{Area}^{-1} \cdot \mbox{time}^{\beta}$, and $[\lambda] =\mbox{time}^{-\alpha}$. Bringing Eq. \eqref{GC-model} to the Laplace domain one gets
\vskip -10pt
	\be
	s \, \wt{J} (s) =
	\frac{a}{b} + \frac{1}{b \, s^\beta (1 - \lambda \, s^{-\alpha})^\gamma}
	= \frac{1}{s \, \wt{G} (s)}
	\ee
\vskip -3pt \noindent
	which, in turn, yields
\vskip -12pt
	\be
	J(t) = \frac{a}{b} + \frac{t^\beta}{b} \, E^{\gamma} _{\alpha, \beta+1} (\lambda \, t^\alpha)  \, ,
	\ee
\vskip -5pt \noindent
	and
\vskip -12pt
	\be
	G(t) = \frac{b}{a} \sum _{n=0} ^{\infty} \left( - \frac{t^\beta}{a} \right)^n \,
	E^{\gamma n} _{\alpha, \beta n + 1} (\lambda \, t^\alpha) \, ,
	\ee
	where the last series is absolutely convergent on $t > 0$ for the whole domain of parameters considered here.
	
	It is also worth mentioning that this model is able to reproduce, as particular cases, the Maxwell, Voigt, and Zener models in both their fractional and ordinary realizations. For details and constraints on the parameter space see \cite{GC-CNSNS2018}. Furthermore, a thorough study of the processes of storage and dissipation of energy in materials described by \eqref{GC-model}
has been carried out by I. Colombaro, A. Giusti, and S. Vitali in \cite{CGV-Mathematics2018}.

\subsection{Stochastic Processes and Diffusion} \label{subsec:stochastic}	

	It is now worth discussing the implications of Prabhakar operators for stochastic processes, and their relation with anomalous diffusion phenomena. For a detailed review on this topic we refer the interested reader to \cite{stanislavsky2018transient,Sandev2019GeneralisedDA,SandevTomovski-Book}.

	It is well-known that the probability density function (PDF) for the location $x$ at time $t$ of a Brownian particle
	satisfies the classical heat equation. In this case, the mean squared displacement (MSD) is linear in time.
	Different alternative non-Brownian models have been studied in the literature in order to describe anomalous diffusions in
	complex media. A widely used and general scheme is given by the Continuous-Time Random Walk (CTRW) approach, first
	introduced in \cite{montroll1965random}. The structure of this stochastic jump process is described by specifying
	the distributions of independent random jump lengths and waiting times (see \eg \cite{klafter2011first}).

	In \cite{Sandev-Diffusion1}, a general class of integro-differential diffusion-type equations given by
\vskip -12pt
	\be
		\label{eq:gen-diff-mod}
		\frac{\partial}{\partial t} W(t, x) =
		\frac{\partial}{\partial t} \int _0^t \eta (t-t') \frac{\partial^2}{\partial{x}^2} W(t' ,x)\, \mathrm{d} t' \, ,
	\ee
	is considered. Here, $W(t,x)$ is the PDF relative to the location of a CTRW admitting jumps of finite variance and waiting times with a
	specific distribution. In particular, the memory kernel $\eta (t)$ in \eqref{eq:gen-diff-mod} is connected to the PDF $\psi(t)$ of the random
	waiting times through their Laplace transforms according to
\vskip -13pt
	\be
		\wt{\psi} (s) = \frac{\wt{\eta} (s)}{1 + \wt{\eta} (s)} \, .
	\ee
	Summarizing, Eq. \eqref{eq:gen-diff-mod} is the governing equation of a family of CTRWs, given some suitable constraints on the memory
	kernel $\eta(t)$.

	One of the interesting cases considered in \cite{Sandev2019GeneralisedDA} is based on the application of the Prabhakar integral.
	In this respect, let us consider, for instance, the case of the {\em tempered time-fractional diffusion equation} with
	$\eta(t)$ such that $\widetilde{\eta}(s) = (s\widetilde{\gamma}(s))^{-1}$, where
\vskip-12pt
	\be
		\label{amode}
		\gamma (t) = \eu ^{-b t} \frac{t^{-\alpha}}{\Gamma(1 - \alpha)} \, ,
	\ee
	with $0 < \alpha < 1$ and $b > 0$ denoting the truncation parameter. The corresponding memory kernel
	then reads
	\be
		\eta(t) = t^{\alpha - 1} E^{\alpha - 1}_{1, \alpha} \big( - b t \big) \, ,
	\ee
	leading to
\vskip -12pt
	\be
		\frac{\partial}{\partial t} W(t, x) =
		\DPrRL ^{1-\alpha} _{1, 1- \alpha, - b;0} \Bigg[ \frac{\partial^2}{\partial{x}^2} W(t,x)\Bigg] \, .
	\ee
	The MSD turns out to be $2 t^{\alpha} E^{\alpha - 1}_{1, \alpha+1} ( - b t)$,
	from which one can immediately infer the emergence of a crossover from subdiffusion to normal diffusion.

	If instead we reverse the problem by considering the case with
	\be
		\eta (t) = \eu ^{-b t} \frac{t^{-\alpha}}{\Gamma(1 - \alpha)},
	\ee
\vskip -3pt \noindent
	then one finds
\vskip -10pt
	\be
		\DPrC ^\alpha _{1, \alpha, -b; 0} W(t, x) =
		\frac{\partial^2}{\partial{x}^2} W(t,x).
	\ee
	Correspondingly, the related MSD becomes
	$2 t^{\alpha} E^{\alpha + 1}_{1, \alpha+1} (-b t)$,
	implying a crossover from subdiffusion to a constant behavior.

	Working in the framework of the theory of CTRWs, it was also shown \cite{Sandev-Diffusion5} that one can build models aimed at characterizing the transition among anomalous diffusion scenarios based on the Prabhakar function and its derivative operators.

	Finally, it is also worth mentioning that a study of the {\em generalized Langevin Equation} with both the Prabhakar and tempered Prabhakar kernels as memory kernels has been carried out by T. Sandev in \cite{Sandev-Diffusion7}.
	
	\section{Applications in renewal processes} \label{sec:probability}

\setcounter{section}{7} \setcounter{equation}{0} 

	Generalizations of renewal processes based on the applications of time-fractional differential operators,
	or special functions related to fractional calculus, have recently attracted much attention.
	Different (non equivalent) approaches to build \textit{fractional} counting processes have been developed in the
	literature. The motivation for the increasing interest on this topic is obviously related to the pervasive applications
	of classical counting processes, like the Poisson process, in real-world models.

	The Mittag--Leffler function appeared as residual waiting time between events in renewal processes already in the 1960s, namely processes with properly scaled thinning out the sequence of events in a power law renewal process (see \cite{gne} and \cite{g3}). A renewal process with Mittag--Leffler--distributed waiting times is in essence a fractional Poisson process.

     Hilfer and Anton \cite{hill2} were the first to introduce the Mittag--Leffler waiting-time density
     \vskip -11pt
\begin{align}
f_{\mu}(t)=-\frac{\d}{\d t}E_{\mu}(-t^{\mu}) = t^{\mu-1} E_{\mu,\mu} (-t^\mu)
\end{align}
in the context of continuous time random walk. Specifically, they were able to show that such a waiting time density is required in order for the evolution of the sojourn density to be governed by a fractional Kolmogorov--Feller equation. In other words, they time-changed a random walk with an independent fractional Poisson process. Systematic studies of both analytic and probabilistic aspects of this topic were then triggered, to the best of our knowledge, in 2000 by the work of O. N. Repin and A. I. Saichev \cite{Repin}. Many (mostly independent) contributions investigated the feature of this new renewal process following the path paved by the previously mentioned pivotal works, see \eg \cite{Beghin,Macci, Laskin, scalas, meer} and references therein. To this regard, it is worth stressing that the Prabhakar function naturally arises in the context of fractional Poisson processes. Indeed, the state probabilities of a time-fractional Poisson process can be expressed in terms of Prabhakar functions, see \eg \cite{Beghin,Laskin} for details.

	Recently, several generalizations of the Poisson process based on the Prabhakar function have appeared in the literature. These extensions clearly offer more flexibility when attempting to capture the main features of real-world renewal processes, also including as special cases the time-fractional and the classical Poisson processes. Essentially, the generalization of the standard Poisson process can be achieved in two different ways: i) deriving the state probabilities of the generalized counting process by solving an infinite system of time-fractional difference-differential equations involving regularized Prabhakar derivatives; ii) considering a renewal process with inter-event time density function involving a Prabhakar function. The first approach was originally developed by R. Garra, {\em et al.} in \cite{GarraGorenfloPolitoTomovski2014}, whereas the second one was first proposed by D. O. Cahoy and F. Polito in \cite{Cahoy2013} (see also \cite{Michelitsch2019, michelitsch2019continuous}). Further, an alternative proposal by T. K. Pog{\'a}ny and {\v{Z}}. Tomovski was provided in \cite{Pogany2016}, where the generalization was obtained in the sense of weighted Poisson distributions.				

	Let us now analyze the first approach, namely the one based on the replacement of the first-order derivative with the regularized Prabhakar derivative. In other words, this procedure has to be implemented onto the infinite system of difference--differential equations governing the state probability of the Poisson process, \ie			
\vskip -12pt
						\begin{align}
							\label{aa1}
							\begin{cases}
								\DPrC^{\gamma}_{\rho, \mu, -\phi} p_k(t) = -\lambda
								p_k(t) +\lambda p_{k-1}(t), & k \ge 0, \: t > 0,
								\: \lambda > 0, \\
								p_k(0) =
								\begin{cases}
									1, \quad k=0, \\
									0, \quad k \geq 1,%
								\end{cases}
								&
							\end{cases}%
						\end{align}
\vskip -4pt \noindent
where $\phi > 0$, $\gamma \ge 0$, $0 < \rho \le 1$, and $0<\mu \le 1$. If $\gamma \ne 0$ one also has that $0 < \mu \lceil \gamma \rceil/\gamma - r\rho < 1$, $\forall \: r=0,\dots,\lceil \gamma \rceil$. These constraints on the parameters are required in order to ensure non-negativity of the solution. Multiplying both terms of \eqref{aa1} by $v^k$, with $v$ denoting an auxiliary variable such that $|v| \le 1$, and summing over all $k$, we obtain the fractional Cauchy problem for the probability generating function $G(v,t) = \sum_{k=0}^\infty v^k p_k(t)$ of the counting number
$N(t)$, $t \ge 0$,
					\begin{equation}
						\label{gen}
						\begin{cases}
						\DPrC^{\gamma}_{\rho, \mu, -\phi} G(v,t)= -\lambda (1-v)G(v,t), &
							|v| \le 1, \\
							G(v,0)=1,  &
						\end{cases}%
					\end{equation}
whose solution reads (see \cite{GarraGorenfloPolitoTomovski2014} for details)
\vskip -10pt
\begin{equation}
\label{G}
G(v,t)=\sum_{k=0}^{\infty}(-\lambda t^{\mu})^k(1-v)^k E_{\rho, \mu k+1}^{\gamma k}(-\phi t^{\rho}),
\qquad |v| \le 1.
\end{equation}
 It is easy to see that, for $\gamma=0$, the latter reduces to
						\begin{equation}
							G(v,t)=\sum_{k=0}^{\infty}\frac{(-\lambda t^{\mu})^k(1-v)^k}{\Gamma(\mu k+1)}
							= E^1_{\mu,1}(-\lambda(1-v)t^{\mu}),
						\end{equation}
that coincides with the probability generating function of the (standard) fractional Poisson process, see \eg \cite{Laskin}.
				
	From the probability generating function \eqref{G} it is possible to infer the probability distribution at fixed time $t \ge 0$
of $N(t)$, governed by the evolution equation \eqref{aa1}. Indeed, a simple binomial expansion leads to
\begin{align}
G(v,t) = \sum_{k =0}^{\infty}v^k \sum_{r=k}^{\infty}(-1)^{r-k}\binom{r}{k}%
(\lambda t^{\mu})^r E_{\rho, \mu r+1}^{\gamma r}(-\phi t^{\rho}) \, ,
\end{align}
\vskip -3pt \noindent
that implies
\vskip -10pt
\begin{equation}
\label{distro}
	p_k(t)=\sum_{r=k}^{\infty}(-1)^{r-k}\binom{r}{k}(\lambda t^{\mu})^r E_{\rho,
			\mu r+1}^{\gamma r}(-\phi t^{\rho}), \qquad k \ge 0, \: t \ge 0.
\end{equation}
From \eqref{aa1}, one can also infer the mean value of $N(t)$. To this aim, it suffices to differentiate Eq. \eqref{gen}
with respect to $v$ and set $v=1$, this leads to
\begin{equation}
\begin{cases}
	\DPrC^{\gamma}_{\rho, \mu, -\phi}\mathbb{E}N(t)=\lambda, & t > 0,
				\\
	\mathbb{E}N(t)\big|_{t=0}=0, &
\end{cases}%
\end{equation}
whose solution is simply given by
\begin{equation}
\label{mean}
\mathbb{E}N(t)=\lambda t^{\mu}E^{\gamma}_{\rho, 1+\mu}(-\phi t^{\rho}),
\qquad t \ge 0 \, .
\end{equation}
It is also possible to prove that the generalized fractional Poisson process $N(t)$ described so far can be constructed as a renewal process with specific waiting times. Indeed, consider $k$ independent and identically distributed (i.i.d.) random variables $T_j$, $j=1,\dots,k$, representing the inter-event waiting times, with probability density functions
\vskip -10pt
\begin{align}
\label{tempi}
f_{T_j}(t_j) = \lambda t_j^{\mu-1} \sum_{r=0}^\infty (-\lambda
t_j^\mu)^r E_{\rho,\mu r+\mu}^{\gamma r+\gamma} (-\phi t_j^\rho),
\qquad t \ge 0, \: \mu \in (0,1) \, .
\end{align}
The probability distribution ${\rm Pr}(N(t)=k)$, $k \geq 0$, can then be obtained by making the renewal structure explicit and comparing the Laplace transform of the state probabilities obtained from the two approaches (see \cite{GarraGorenfloPolitoTomovski2014} for details).

Lastly, it is worth remarking that the generalized fractional Poisson process $N(t)$, similarly to the classical fractional Poisson process, can be represented as a time-changed Poisson process in which the random change of time is given by a suitable independent stochastic process with increasing paths. See Section 5.3.1 of \cite{GarraGorenfloPolitoTomovski2014} for details.
					
	The second approach, based on the construction of a renewal process with inter-event time density function involving a Prabhakar function, was originally proposed in \cite{Cahoy2013}.  In this case the generalized Poisson process $N^{\nu,\delta}(t)$ denotes a renewal process with i.i.d. random waiting times distributed as
\vskip -10pt
\begin{equation}
    f^{\nu,\delta}(t) = \lambda^\delta t^{\nu\delta-1}E^{\delta}_{\nu,\delta\nu}(-\lambda t^\nu), \quad \lambda>0,\quad \delta\in \mathbb{R}, \ \nu\in (0,1].
\end{equation}
\vskip -2pt \noindent
Clearly, the latter is a generalization of the (canonical) time-fractional Poisson process, recovered for $\delta = 1$, and of the classical Poisson process, which is instead reproduced when setting $\delta = \nu = 1$. Then, by means of the Laplace transform method, it is possible to prove that the state probabilities $p^{\nu,\delta}_k(t)={\rm Pr} (N^{\nu,\delta}(t)=k)$ read
					    \begin{equation}
					    p_k^{\nu,\delta}(t)= \lambda^{\delta k} t^{\nu\delta k}E^{\delta k}_{\nu,\nu\delta k+1}(-\lambda t^\nu)-\lambda^{\delta (k+1)} t^{\nu\delta (k+1)}E^{\delta (k+1)}_{\nu,\nu\delta(k+1)+1}(-\lambda t^\nu)
					    \end{equation}
and satisfy the Volterra equation
\vskip -10pt
\begin{equation}
p_k^{\nu,\delta}(t) = \lambda^\delta \left(\JPr^\delta_{\nu,\nu\delta,-\lambda;0^+} p_{k-1}^{\nu,\delta}\right)(t)
\end{equation}
\vskip -2pt \noindent
involving the Prabhakar fractional integral.
					
	Finally, we briefly recall that a further proposal for the construction of a generalized Prabhakar-Poisson distribution, based on the Poisson distribution approach, was developed in \cite{Pogany2016}. The probability mass function of a weighted
Poisson process is of the form (we adopt the notation in \cite{Bala2008} and references therein)
\vskip -12pt
\begin{equation}
P\{N^w(t)=n\}= \frac{w(n)p(n,x)}{\mathbb{E}[w(N)]},\quad n\geq 0,
\end{equation}
\vskip -3pt \noindent
where $N$ is a random variable with a Poisson distribution 
$p(n,x)= (x^n/n!)$ $ \eu ^{-x}$, $w(\cdot)$ is a non-negative weight function with non-zero finite expected value, \ie
\vskip -12pt
\begin{equation}
 0<\mathbb{E}[w(N)]= \sum_{n=0}^{\infty}w(n) p(n,x) <\infty.
\end{equation}
Therefore, it is possible to build a Prabhakar-based generalization of the Poisson distribution by choosing the weights as
$w(n,\gamma,\alpha,\beta)= \Gamma(n+\gamma)/[\Gamma(\gamma)\Gamma(\alpha n+\beta)]$. This procedure provides a generalization of the Poisson distribution, but it does not represent a renewal process unlike the other approaches. On the other hand, this last generalization can be useful in the context of sub- or super-Poissonian distribution applications.

\vspace*{-3pt}

\section{Numerical aspects} \label{sec:numerics}

\setcounter{section}{8} \setcounter{equation}{0} 

	The numerical evaluation of the Prabhakar function is a quite delicate and involved topic that appears to have been loosely
discussed in the current literature, except in a couple of works in \cite{Garrappa2015_SIAM,StanislavskyWeron2012}.

The simplest way to numerically evaluate the Prabhakar function consists in relaying on its definition \eqref{eq:PrabhakarFunction}. Namely, for such an approach it would be sufficient to fix a sufficiently large, though finite, number $K\in\mathbb{N}$ and approximate the function by the truncated series
\vskip -10pt
\begin{equation}\label{eq:ApproxSeries}
	E_{\alpha,\beta}^{\gamma}(z) \approx \frac{1}{\Gamma(\gamma)}\sum_{k=0}^{K} \frac{ \Gamma(\gamma+k) z^{k}}{k! \Gamma(\alpha k + \beta)}.
\end{equation}
\vskip -3pt \noindent
However, some critical issues emerge when this simple scheme is implemented. Indeed, the Gamma function grows very fast as its argument gets larger. This means, from a numerical perspective, that a restriction on $K$ has to be imposed to avoid overflow, \ie the generation of numbers exceeding the range that can be represented with a given number of digits. In the standard IEEE-754 double precision arithmetic of a commonly used computer $1.8 \times 10^{308}$ is approximately the largest floating-point number that can be represented. Thus, since $\Gamma(171.624)\approx 1.8 \times 10^{308}$, the maximum number of terms in the truncated series \eqref{eq:ApproxSeries} is bounded by $K < {(171.624 - \beta)}/{\alpha}$. This bound heavily restricts the applicability of this approach since the series \eqref{eq:PrabhakarFunction} converges very slowly for arguments with moderate or large modulus. Furthermore, additional problems arise when trying to deal with arguments outside the positive real axis, in particular when $|z|>1$. In this case, as the value of $k$ increases the sum \eqref{eq:ApproxSeries} has consecutive terms with large modulus but opposite sign whose sum is extremely ill-conditioned in the finite-precision arithmetic of computers. As a consequence, one often incurs in catastrophic numerical cancellation. Summing up, the approximation \eqref{eq:ApproxSeries} can be applied only for arguments with small, or just very moderately large, modulus. Otherwise, this approach turns out to be unreliable. This is, after all, the same conclusion one reaches by carrying out a similar analysis for the Mittag--Leffer function and its derivatives (see \eg \cite{GorenfloLoutchkoLuchko2002, GarrappaPopolizio2018}).

	A viable strategy for the numerical evaluation of the Prabhakar function for arguments with large modulus is to relay on some of the asymptotic expansions presented in Section \ref{SS:AsymptoticBehavior}. Whereas, the  remaining cases can be addressed by working backward from the Laplace domain.
	
Regarding the method based on the Laplace transform, the very simple expression given in (\ref{ML3_LT}) suggests to numerically compute the Prabhakar function by means of the inversion of its Laplace transform. J. A. C. Weideman and L. N. Trefethen paved the way for the current developments of this approach, analyzing the case $\beta=\gamma=1$ \cite{WeidemanTrefethen2007}. Later on, the analysis was extended to any $\beta$ and $\gamma=1$ \cite{GarrappaPopolizio2013} and, successively, to any $\gamma>0$ for $|\arg z|> \alpha \pi$ and $0<\alpha<1$ \cite{Garrappa2015_SIAM}. Then, these results were further extended to more general arguments $z$, with $\gamma\in\mathbb{N}$ in \cite{GarrappaPopolizio2018}.

To give just a general idea of this last technique, the starting point is represented by the expression
\vskip -10pt
\begin{equation}\label{eq:Pb_InvLT}
E_{\alpha,\beta}^{\gamma}(z) = \frac{1}{2\pi i} \int_{\bar{\mathcal{C}}} e^{s} H(s;z) \, \du s
, \quad H(s;z) = \frac{s^{\alpha\gamma - \beta}}{(s^{\alpha} - z)^{\gamma}} ,
\end{equation}
obtained from the formula for the inversion of the Laplace transform, where the considered contour $\bar{\mathcal{C}}$ is a parabola that begins and ends in
the left half-plane, such that $\real(s)\rightarrow -\infty$ at each end. In this way, the integrand in (\ref{eq:Pb_InvLT}) rapidly decays as  $\real(s)\rightarrow -\infty$, being forced by the exponential factor. Moreover, the contour $\bar{\mathcal{C}}$ must encompass all the singularities of $H(s;z)$ (that would all be lying on the left-hand side of the vertex of $\bar{\mathcal{C}}$) and the branch-cut imposed on the negative real semi axis to make $H(s;z)$ single-valued. All that being said, the trapezoidal quadrature rule can now be profitably applied for the numerical evaluation of \eqref{eq:Pb_InvLT}. However, this procedure requires to accurately fix some parameters. Indeed, in order to define the contour $\bar{\mathcal{C}}$ one parameter needs to be set, while the trapezoidal rule requires fixing the number of nodes and the spacing between them. In practice, if
\vskip -10pt
\[
w(u) = \mu (iu+1)^2,\,\, -\infty < u < \infty,
\]
describes the parabolic contour, $h>0$ is a given step-size and the quadrature nodes are defined as $u_k = kh$ for $k=-N,\ldots,N$, then the resulting approximation to \eqref{eq:Pb_InvLT} reads
\vskip -11pt
\begin{equation}\label{eq:approxInvLT}
{\tilde{E}}^{[N,h]}(t) = \frac{h}{2 \pi i} \sum_{k=-N}^{N}  e^{w(u_k)} H (w(u_k);t) w^{\prime}(u_k) .
\end{equation}
To choose the optimal parameters $N$, $h$ and $\mu$ a deep error analysis must be carried out by requiring that the global error is proportional to the required accuracy (usually something close to the machine precision to force the method to provide the highest possible accuracy). Specifically, the influence of the parameters $N$, $h$ and $\mu$ on the round-off, truncation, and discretization errors have to be considered. Basically, this means that one needs to find the most suitable contour $\bar{\mathcal{C}}$. Incidentally, such a contour, roughly speaking, turns out to be the one whose branches are the farthest from the singularities of the integrand. Now, once the required accuracy is fixed, the algorithm finds the optimal parameters and then computes \eqref{eq:approxInvLT}. In particular, the value $N$ is kept as small as possible and, since it determines the computational cost of the procedure, the resulting algorithm turns out to be fast and very accurate (see \cite{Garrappa2015_SIAM} for details).
		
If the singularities of $H(s;t)$ turn out to be to scattered in the complex plane and too far away from both the origin and the branch-cut imposed on the negative real semi axis, it is essential to allow the contour to go over some of the singularities and apply to \eqref{eq:Pb_InvLT} the residue subtraction formula
\[
	E^\gamma _{\alpha,\beta}(z) =
		\sum_{s^{\star} \in S^{\star}_{\bar{\mathcal{C}}}} \operatorname{Res} \bigl( e^{s} H(s;z) , s^{\star} \bigr) +
		\frac{1}{2\pi i} \int_{\bar{\mathcal{C}}} e^{s} H(s;z)  \du s,
\]
with $S^{\star}_{\bar{\mathcal{C}}}$ the set of the singularities of $H(s;z)$ laying on the right-hand side of the contour $\bar{\mathcal{C}}$. This correction can, however, be implemented just for $\gamma\in\mathbb{N}$, otherwise the more involved branch-cut due to the two real powers in $H(s,z)$ does not allow the contour to pass over any of the singularities. This is the reason for which the method described in \cite{Garrappa2015_SIAM}, and corresponding Matlab code \cite{Garrappa2014ml}, are restricted to $0<\alpha<1$ and $|\arg z| > \alpha \pi$, whilst in \cite{GarrappaPopolizio2018} it is discussed the case $\gamma\in\mathbb{N}$ and a formula for the residues $\operatorname{Res} \bigl( e^{s} H(s;z) , s^{\star} \bigr)$ is given.

 		
	\section{Discussion and outlook} \label{sec:conclusions}

\setcounter{section}{9} \setcounter{equation}{0} 

	And here is where this Hitchhiker's Guide to the wonders of Prabhakar fractional calculus comes to an end. 
In Section \ref{sec:history} we have offered a brief recap of the main events that have led to the discovery of the Prabhakar function and the associated fractional operators. In Section \ref{sec:motivations} we have shown how this generalization was not just a mathematical curiosity, but rather something that physics was begging for in order to explain some peculiar phenomena. After that, in Section \ref{sec:math} we have provided an extensive description of the mathematical features of the Prabhakar function. Then, in Section \ref{sec:pracalculus} we combined the results in Section \ref{sec:math} with the general wisdom of the standard (Riemann--Liouville and Caputo) approaches to fractional calculus in order to reconstruct Prabhakar's theory. In Section \ref{sec:physics} we have then framed the physical problems raised in Section \ref{sec:motivations} within this new scheme, thus summarizing the main physical application of this general set-up. In Section \ref{sec:probability} we summarized how Prabhakar's theory can be applied in probability theory, with particular regard for the theory of renewal processes. Finally, in Section \ref{sec:numerics} we have discussed the state of the art of the research on the numerical evaluation of the Prabhakar function.

	From a physical perspective, setting aside the many possible scenarios in which Prabhakar's theory might potentially emerge as a powerful modeling tool, it is worth stressing that much of the fundamental physics connected with the phenomena discussed in Section \ref{sec:physics} is still uncharted territory. First, while this mathematical framework seems to be the key for the mathematical description of anomalous dielectrics of the Havriliak--Negami type, a first-principle derivation of these physical behaviors is still missing. Second, even though the study of viscoelastic applications of Prabhakar calculus, initiated by A. Giusti and I. Colombaro with \cite{GC-CNSNS2018}, has attracted some attention in the literature, most physical aspects of this extension of linear viscoelasticity are yet to be unveiled. These two research lines represent the main open questions concerning the applications of Prabhakar fractional calculus to ``real-world problems''.

	From a probabilistic point of view, we can recall that in \cite{Michelitsch2019, michelitsch2019continuous} the authors analyzed generalized fractional diffusion models on undirected graphs and infinite lattices. This application, making use of Prabhakar random variables and therefore of Prabhakar functions, seems to be a promising topic of research in the context of non-Markovian models on networks. Moreover, further investigations on applications of Prabhakar integrals and derivatives to more general or different point processes, such as birth-death processes or point processes related to macroevolutionary models (see \eg \cite{Modern-stochastics-Polito}), are surely worth pursuing. Finally, according to \cite{KonnoTamura2018}, the applications of fractional renewal processes based on Prabhakar functions to models for neural spiking events could represent compelling research subjects.
	
	From a numerical viewpoint, a further boost should be given toward the development of methods for
the computation of the Prabhakar function, in order to cover a larger range of parameters. Moreover,
another stimulating research topic, that supposedly will attract much attention in the future, concerns
the possibility of developing highly efficient methods for solving differential equations with Prabhakar derivatives,
particularly for applications in computational electromagnetism (see \eg \cite{BiaCaratelliMesciaCicchettiMaionePrudenzano2014,
CausleyPetropoulosJiang2011,GarrappaMaione2017}).

\vspace*{-2pt}

\section*{Acknowledgments}	

A. Giusti is supported by the Natural Sciences and Engineering Research Council of Canada (Grant No.~2016-03803 to V.F.) and by Bishop's University. 
I. Colombaro is supported by the Spanish Ministry of Science, Innovation and Universities under grants TEC2016-78434-C3-1-R and BES-2017-081360. 
R. Garrappa is partially supported by the COST Action CA 15225 ``Fractional-order systems-analysis, synthesis and their importance for future design'' and by an INdAM-GNCS 2019 project. 
F. Polito has been partially supported by the project ``Memory in Evolving Graphs'' (Compagnia di San Paolo--Universit{\`a} di Torino). 
The work of M. Popolizio is supported in part by an INdAM-GNCS 2019 project.

This work has been carried out in the framework of the activities of the Italian National Group of Mathematical Physics (INdAM-GNFM) and of the Italian National Group for Scientific Computing (INdAM-GNCS).
	
\vspace*{-2pt}

\section*{Authors' Contributions}

\textbf{Andrea Giusti} has conceived the idea, organized, and coordinated the work. He has collected the contributions of the authors and he actively worked on Sections \ref{sec:motivations}, \ref{sec:math}, \ref{sec:pracalculus}, and \ref{sec:physics}. He has also prepared the final version of the text. \textbf{Ivano Colombaro} has contributed to the content of Sections \ref{sec:motivations} and \ref{sec:physics}. \textbf{Roberto Garra} has coordinated and contributed to the work on Sections \ref{sec:math} and \ref{sec:probability}. He also offered valuable insights on the content of Section \ref{subsec:stochastic}. \textbf{Roberto Garrappa} has coordinated and worked on Sections \ref{sec:numerics}. He has also contributed to the work on Sections \ref{sec:history}, \ref{sec:math}, and \ref{sec:pracalculus}. \textbf{Federico Polito} has contributed to the work on Sections \ref{sec:math}, \ref{sec:probability}, and he has offered valuable insights on the content of Section \ref{subsec:stochastic}. \textbf{Marina Popolizio} has worked on Section \ref{sec:numerics}. \textbf{Francesco Mainardi} has contributed to Section \ref{sec:history} and he has offered his guidance and insight to all the parts appearing in this work. Finally, all authors have approved the content of the manuscript.



\begin{thebibliography}{99} 
\normalsize 


\bibitem{Bala2008}
N.~Balakrishnan, T.~J. Kozubowski, {A class of weighted Poisson processes}.
  \emph{Stat. Probabil. Lett.} \textbf{78}, No~15 (2008), 2346--2352.

\bibitem{Barrett1954}
J.~H. Barrett, {Differential equations of non-integer order}. \emph{Canadian J.
  Math.} \textbf{6} (1954), 529--541.

\bibitem{Macci}
L.~Beghin, C.~Macci, {Large deviations for fractional Poisson processes}.
  \emph{Stat. Probabil. Lett.} \textbf{83}, No~4 (2012), 1193--1202.

\bibitem{Beghin}
L.~Beghin, E.~Orsingher, {Fractional Poisson processes and related planar
  random motions}. \emph{Electron. J. Probab.} \textbf{14}, No~61 (2009),
  1790--1826.

\bibitem{BeghinOrsingher2010}
L.~Beghin, E.~Orsingher, {Poisson-type processes governed by fractional and
  higher-order recursive differential equations}. \emph{Electron. J. Probab.}
  \textbf{15}, No~22 (2010), 684--709.

\bibitem{BiaCaratelliMesciaCicchettiMaionePrudenzano2014}
P.~Bia, D.~Caratelli, L.~Mescia, R.~Cicchetti, G.~Maione, F.~Prudenzano, {A
  novel {FDTD} formulation based on fractional derivatives for dispersive
  {H}avriliak--{N}egami media}. \emph{Signal Process.} \textbf{107} (2015),
  312--318.

\bibitem{Braaksma1964}
B.~L.~J. Braaksma, {Asymptotic expansions and analytic continuations for a
  class of {B}arnes-integrals}. \emph{Compositio Math.} \textbf{15} (1964),
  239--341.

\bibitem{Buhl1925}
A.~Buhl, \emph{{S\'eries analytiques. Sommabilit\'e}}. Number~7 in M\'emorial
  des sciences math\'ematiques, Gauthier-Villars (1925).

\bibitem{Cahoy2013}
D.~O. Cahoy, F.~Polito, {Renewal processes based on generalized Mittag-Leffler
  waiting times}. \emph{Commun. Nonlin. Sci. Numer. Simul.} \textbf{18}, No~3 (2013), 639--650.

\bibitem{CamargoCharnetCapelasDeOliveira2009}
R.~F. Camargo, R.~Charnet, E.~Capelas~de Oliveira, {On some fractional Green's
  functions}. \emph{J. Math. Phys.} \textbf{50}, No~4 (2009), 043514.

\bibitem{CamargoChiacchioCharnetCapelasDeOliveira2009}
R.~F. Camargo, A.~O. Chiacchio, R.~Charnet, E.~C. de~Oliveira, {Solution of the
  fractional {L}angevin equation and the {M}ittag-{L}effler functions}.
  \emph{J. Math. Phys.} \textbf{50}, No~6 (2009), 063507.

\bibitem{CapelasMainardiVaz2011}
E.~Capelas De~Oliveira, F.~Mainardi, J.~Vaz~Jr, {Models based on
  {M}ittag-{L}effler functions for anomalous relaxation in dielectrics}.
  \emph{Eur. Phys. J. Spec. Top.} \textbf{193}, No~1 (2011), 161--171 [Revised
  version: \href{https://arxiv.org/abs/1106.1761v2}{arXiv:1106.1761v2}].

\bibitem{CaputoMainardi1971-2}
M.~Caputo, F.~Mainardi, {A new dissipation model based on memory mechanism}.
  \emph{Pure Appl. Geophys.} \textbf{91}, No~1 (1971), 134--147.

\bibitem{CaputoMainardi1971-1}
M.~Caputo, F.~Mainardi, {Linear models of dissipation in anelastic solids}.
  \emph{Riv. del Nuovo Cim.} \textbf{1}, No~2 (1971), 161--198.

\bibitem{CausleyPetropoulosJiang2011}
M.~F. Causley, P.~G. Petropoulos, S.~Jiang, {Incorporating the
  {H}avriliak-{N}egami dielectric model in the {FD}-{TD} method}. 
  \emph{J. Comput. Phys.} \textbf{230}, No~10 (2011), 3884--3899.

\bibitem{ChamatiTonchev2006}
H.~Chamati, N.~S. Tonchev, {Generalized {M}ittag-{L}effler functions in the
  theory of finite-size scaling for systems with strong anisotropy and/or
  long-range interaction}. \emph{J. Phys. A: Math. Gen.} \textbf{39}, No~3
  (2005), 469--478.

\bibitem{Cole1933electric}
K.~S. Cole, {Electric conductance of biological systems}. In: \emph{Cold Spring
  Harbor Symposia on Quantitative Biology}, Vol.~1, 107--116, Cold Spring
  Harbor Laboratory Press (1933).

\bibitem{CC1}
K.~S. Cole, R.~H. Cole, {Dispersion and Absorption in Dielectrics I.
  Alternating Current Characteristics}. \emph{J. Chem. Phys.} \textbf{9}
  (1941), 341.

\bibitem{CC2}
K.~S. Cole, R.~H. Cole, {Dispersion and Absorption in Dielectrics II. Direct
  Current Characteristics}. \emph{J. Chem. Phys.} \textbf{10} (1942), 98.

\bibitem{Bessel-2017}
I.~Colombaro, A.~Giusti, F.~Mainardi, {A class of linear viscoelastic models
  based on Bessel functions}. \emph{Meccanica} \textbf{52}, No 4-5 (2017), 825--832.

\bibitem{Bessel-2017-2}
I.~Colombaro, A.~Giusti, F.~Mainardi, {On the propagation of transient waves in
  a viscoelastic Bessel medium}. \emph{Z. Angew. Math. Phys.} \textbf{68}, No~3
  (2017), \# 62.

\bibitem{WaveMotion-2017}
I.~Colombaro, A.~Giusti, F.~Mainardi, {On Transient Waves in Linear
  Viscoelasticity}. \emph{Wave Motion} \textbf{74} (2017), 191--212.

\bibitem{CGV-Mathematics2018}
I.~Colombaro, A.~Giusti, S.~Vitali, {Storage and dissipation of energy in
  Prabhakar viscoelasticity}. \emph{Mathematics} \textbf{6} (2018), 15.

\bibitem{DC}
D.~W. Davidson, R.~H. Cole, {Dielectric relaxation in glycerol, propylene
  glycol and n-propanol}. \emph{J. Chem. Phys.} \textbf{19} (1951), 1484--1491.

\bibitem{Davis1936theory}
H.~T. Davis, \emph{{Handbook of Mathematical Functions: With Formulas, Graphs,
  and Mathematical Tables}}. Bloomington, Ind. (1936).

\bibitem{Debye-Seminal}
P.~Debye, {Zur theorie der spezifischen W\''{a}rme}. \emph{Ann. Phys.}
  \textbf{39} (1912), 789--839.

\bibitem{Diethelm2010}
K.~Diethelm, \emph{{The Analysis of Fractional Differential Equations}}. Vol.
  2004 of \emph{Lecture Notes in Mathematics}, Springer-Verlag, Berlin (2010).

\bibitem{Djrbashian1993}
M.~M. Djrbashian, \emph{{Harmonic Analysis and Boundary Value Problems in the
  Complex Domain}}. Vol.~65 of \emph{Operator Theory: Advances and
  Applications}, Birkh\"{a}user Verlag, Basel (1993) [Transl. from the
  manuscript by H.M. Jerbashian and A.M. Jerbashian (A.M. Dzhrbashyan)].

\bibitem{DOvidioPolito2018}
M.~D'Ovidio, F.~Polito, {Fractional diffusion--telegraph equations and their
  associated stochastic solutions}. \emph{Theory Probab. Appl.} \textbf{62},
  No~4 (2018), 552--574 [arXiv: 1307.1696 (2013)].

\bibitem{ErdelyiMagnusOberhettingerTricomi3}
A.~Erd\'{e}lyi, W.~Magnus, F.~Oberhettinger, F.~G. Tricomi, \emph{Higher
  Transcendental Functions. {V}ol. {III}}. McGraw-Hill Book Company, Inc., New
  York-Toronto-London (1955). Based, in part, on notes left by Harry Bateman.

\bibitem{Fox1928}
C.~Fox, {The asymptotic expansion of integral functions defined by generalized
  hypergeometric functionss}. \emph{Proc. London Math. Soc.} \textbf{s2--27},
  No~1 (1928), 389--400.

\bibitem{GarraGarrappa2018}
R.~Garra, R.~Garrappa, The {P}rabhakar or three parameter {M}ittag-{L}effler
  function: Theory and application. \emph{Commun. Nonlin. Sci. Numer. Simul.}
  \textbf{56} (2018), 314--329.

\bibitem{GarraGorenfloPolitoTomovski2014}
R.~Garra, R.~Gorenflo, F.~Polito, {\v{Z}}.~Tomovski, Hilfer-{P}rabhakar
  derivatives and some applications. \emph{Appl. Math. Comput.} \textbf{242}
  (2014), 576--589.

\bibitem{Garrappa2014ml}
R.~Garrappa, The {M}ittag--{L}effler function. \emph{MATLAB Central File Exchange}
  (2014), File ID: 48154.

\bibitem{Garrappa2015_SIAM}
R.~Garrappa, Numerical evaluation of two and three parameter {M}ittag-{L}effler
  functions. \emph{SIAM J. Numer. Anal.} \textbf{53}, No~3 (2015), 1350--1369.

\bibitem{GarrappaHN2016}
R.~Garrappa, Gr\"{u}nwald-{L}etnikov operators for fractional relaxation in
  {H}avriliak-{N}egami models. \emph{Commun. Nonlin. Sci. Numer. Simul.}
  \textbf{38} (2016), 178--191.

\bibitem{GarrappaMainardiMaione2016}
R.~Garrappa, F.~Mainardi, G.~Maione, Models of dielectric relaxation based on
  completely monotone functions. \emph{Fract. Calc. Appl. Anal.} \textbf{19},
  No~5 (2016), 1105--1160.

\bibitem{GarrappaMaione2017}
R.~Garrappa, G.~Maione, {Fractional {P}rabhakar derivative and applications in
  anomalous dielectrics: A numerical approach}. \emph{Lecture Notes in
  Electrical Engineering} \textbf{407} (2017), 429--439.

\bibitem{GarrappaPopolizio2013}
R.~Garrappa, M.~Popolizio, Evaluation of generalized {M}ittag--{L}effler
  functions on the real line. \emph{Adv. Comput. Math.} \textbf{39}, No~1
  (2013), 205--225.

\bibitem{GarrappaPopolizio2018}
R.~Garrappa, M.~Popolizio, {Computing the matrix Mittag--Leffler function with
  applications to fractional calculus}. \emph{J. Sci. Comput.} \textbf{77},
  No~1 (2018), 129--153.

\bibitem{Giusti2018_NODY}
A.~Giusti, A comment on some new definitions of fractional derivative.
  \emph{Nonlinear Dyn.} \textbf{93}, No~3 (2018), 1757--1763.

\bibitem{Giusti2019}
A.~Giusti, {General fractional calculus and Prabhakar's theory}. \emph{Commun.
  Nonlin. Sci. Numer. Simul.} \textbf{83} (2020), 105114.

\bibitem{GC-CNSNS2018}
A.~Giusti, I.~Colombaro, {Prabhakar-like fractional viscoelasticity}.
  \emph{Commun. Nonlin. Sci. Numer. Simul.} \textbf{56} (2018), 138--143.

\bibitem{GiustiMainardi-2016}
A.~Giusti, F.~Mainardi, {A dynamic viscoelastic analogy for fluid-filled
  elastic tubes}. \emph{Meccanica} \textbf{51}, No~10 (2016), 2321--2330.

\bibitem{Gross-Fuoss-GM}
A.~Giusti, F.~Mainardi, {On infinite series concerning zeros of Bessel
  functions of the first kind}. \emph{Eur. Phys. J. Plus} \textbf{131}, No~6
  (2016), \# 206.

\bibitem{gne}
B.~V. Gnedenko, I.~N. Kovalenko, \emph{Introduction to Queueing Theory}. Israel
  Program for Scientific Translations, Jerusalem (1968).

\bibitem{GorenfloKilbasMainardiRogosin2014}
R.~Gorenflo, A.~A. Kilbas, F.~Mainardi, S.~Rogosin, \emph{Mittag-Leffler
  functions. Theory and Applications}. Springer Monographs in Mathematics,
  Springer, Berlin (2014).

\bibitem{GorenfloLoutchkoLuchko2002}
R.~Gorenflo, J.~Loutchko, Y.~Luchko, Computation of the {M}ittag-{L}effler
  function {$E_{\alpha,\beta}(z)$} and its derivative. \emph{Fract. Calc. Appl.
  Anal.} \textbf{5}, No~4 (2002), 491--518.

\bibitem{GorskaPLA2019}
K.~G\'{o}rska, A.~Horzela, L.~A, {Composition law for the Cole-Cole relaxation
  and ensuing evolution equations}. \emph{Phys. Lett. A} \textbf{383}, No~15
  (2019), 1716--1721.

\bibitem{GorskaJPA2018}
K.~G\'{o}rska, A.~Horzela, G.~Dattoli, P.~K. A, {The Havriliak-Negami
  relaxation and its relatives: the response, relaxation and probability
  density functions}. \emph{J. Phys. A} \textbf{51}, No~13 (2018), 135202.

\bibitem{GorskaHorzelaLattanziPogani2018}
K.~G\'{o}rska, A.~Horzela, T.~K. Lattanzi, A.~Pog\'{a}ny, On the complete
  monotonicity of the three parameter generalized {M}ittag-{L}effler function
  $e_{\alpha,\beta}^{\gamma}(-x)$. Available as arXiv: 1811.10441 (2018).

\bibitem{Gorska-etc-FCAA2019} 
K.~G\'{o}rska, A.~Horzela,  R.~Garrappa,
Some results on the complete monotonicity of Mittag-Leffler functions of Le Roy type.
\emph{Fract. Calc. Appl. Anal.} {\bf 22}, No~5 (2019), 1284--1306.

\bibitem{GorskaKochubei}
K.~G\'{o}rska, A.~Horzela, T.~K. Pog\'{a}ny, {A note on the article ``Anomalous
  relaxation model based on the fractional derivative with a Prabhakar-like
  kernel''}. \emph{Z. Angew. Math. Phys.} \textbf{70} (2019), 141.

\bibitem{Gross1947creep1}
B.~Gross, {On creep and relaxation}. \emph{J. Appl. Phys.} \textbf{18}, No~2
  (1947), 212--221.

\bibitem{Gross1947creep2}
B.~Gross, {On Creep and Relaxation, II.} \emph{J. Appl. Phys.} \textbf{19},
  No~3 (1948), 257--264.

\bibitem{Gross-Fuoss-3}
B.~Gross, {Electrical analogs for viscoelastic systems}. \emph{J. Polym. Sci.}
  \textbf{20}, No~95 (1956), 371--380.

\bibitem{Gross-Fuoss-2}
B.~Gross, {Ladder structures for representation of viscoelastic systems, II}.
  \emph{J. Polym. Sci.} \textbf{20}, No~94 (1956), 123--131.

\bibitem{Gross-Fuoss-1}
B.~Gross, R.~M. Fuoss, {Ladder structures for representation of viscoelastic
  systems}. \emph{J. Polym. Sci.} \textbf{19}, No~91 (1956), 39--50.

\bibitem{Gross1947creep3}
B.~Gross, H.~Pelzer, {On creep and relaxation, III}. \emph{J. Appl. Phys.}
  \textbf{22}, No~8 (1951), 1035--1039.

\bibitem{GurtinSternberg1962}
M.~E. Gurtin, E.~Sternberg, {On the linear theory of viscoelasticity}.
  \emph{Arch. Ration. Mech. Anal.} \textbf{11}, No~1 (1962), 291--356.

\bibitem{HanygaSeredynska2008}
A.~Hanyga, M.~Seredy\'{n}ska, {On a mathematical framework for the constitutive
  equations of anisotropic dielectric relaxation}. \emph{J. Stat. Phys.}
  \textbf{131} (2008), 269--303.

\bibitem{HauboldMathaiSaxena2011}
H.~J. Haubold, A.~M. Mathai, R.~K. Saxena, Mittag-Leffler functions and their
  applications. \emph{J. Appl. Math.} \textbf{2011} (2011), 298628. 

\bibitem{HN}
S.~Havriliak, S.~Negami, {A complex plane analysis of $\alpha$-dispersions in
  some polymer systems}. \emph{J. Polym. Sci. C} \textbf{14} (1966), 99--117.

\bibitem{HavriliakNegami1967}
S.~Havriliak, S.~Negami, A complex plane representation of dielectric and
  mechanical relaxation processes in some polymers. \emph{Polymer} \textbf{8}
  (1967), 161--210.

\bibitem{HavriliakNegami1969a}
S.~Havriliak, S.~Negami, On the equivalence of dielectric and mechanical
  dispersions in poly($n$-hexyl methacrylate). \emph{J. Phys. D Appl. Phys.}
  \textbf{2}, No~9 (1969), 1301--1315.

\bibitem{HavriliakNegami1969b}
S.~Havriliak, S.~Negami, On the equivalence of dielectric and mechanical
  dispersions in some polymers; e.g. poly($n$-octyl methacrylate).
  \emph{Polymer} \textbf{10} (1969), 859 -- 872.

\bibitem{HH1994}
S.~Havriliak~Jr, S.~J. Havriliak, {Results from an unbiased analysis of nearly
  1000 sets of relaxation data}. \emph{J. Non-Cryst. Solids} \textbf{172--174,
  Part 1} (1994), 297--310.

\bibitem{Hilfer2002b}
R.~Hilfer, {Experimental evidence for fractional time evolution in glass
  forming materials}. \emph{Chem. Phys.} \textbf{284}, No~1 (2002), 399--408.

\bibitem{Hilfer2002a}
R.~Hilfer, {H}-function representations for stretched exponential relaxation
  and non-{D}ebye susceptibilities in glassy systems. \emph{Phys. Rev. E}
  \textbf{65} (2002), 061510.

\bibitem{hill2}
R.~Hilfer, L.~Anton, {Fractional master equation and fractal time random
  walks}. \emph{Phys. Rev. E} \textbf{51} (1995), R848--R851.

\bibitem{Debye}
R.~M. Hill, L.~A. Dissado, {Debye and non-Debye relaxation}. \emph{J. Phys. C}
  \textbf{18}, No~19 (1985), 3829.

\bibitem{HilleTamarkin1930}
E.~Hille, J.~D. Tamarkin, On the theory of linear integral equations.
  \emph{Ann. of Math. (2)} \textbf{31}, No~3 (1930), 479--528.

\bibitem{Humbert1953}
P.~Humbert, Quelques r\'{e}sultats relatifs \`a la fonction de
  {M}ittag-{L}effler. \emph{C. R. Acad. Sci. Paris} \textbf{236} (1953),
  1467--1468.

\bibitem{HumbertAgarwal1953}
P.~Humbert, R.~P. Agarwal, Sur la fonction de {M}ittag-{L}effler et
  quelques-unes de ses g\'{e}n\'{e}ralisations. \emph{Bull. Sci. Math. (2)}
  \textbf{77} (1953), 180--185.

\bibitem{HumbertDelerue1953}
P.~Humbert, P.~Delerue, Sur une extension \`a deux variables de la fonction de
  {M}ittag-{L}effler. \emph{C. R. Acad. Sci. Paris} \textbf{237} (1953),
  1059--1060.

\bibitem{Jackson}
J.~D. Jackson, \emph{{Classical Electrodynamics}}. John Wiley \& Sons Inc.
  (1998).

\bibitem{exp1}
A.~K. Jonscher, {The universal dielectric response}. \emph{Nature}
  \textbf{267}, No 5613 (1977), 673--679.

\bibitem{Jonscher1999}
A.~K. Jonscher, {Dielectric relaxation in solids}. \emph{J. Phys. D}
  \textbf{32}, No~14 (1999), R57--R70.

\bibitem{KilbasSaigoSaxena2002}
A.~A. Kilbas, M.~Saigo, R.~K. Saxena, Solution of {V}olterra
  integrodifferential equations with generalized {M}ittag-{L}effler function in
  the kernels. \emph{J. Integral Equations Appl.} \textbf{14}, No~4 (2002),
  377--396.

\bibitem{KilbasSaigoSaxena2004}
A.~A. Kilbas, M.~Saigo, R.~K. Saxena, Generalized {M}ittag-{L}effler function
  and generalized fractional calculus operators. \emph{Integr. Transf.
  Spec. Funct.} \textbf{15}, No~1 (2004), 31--49.

\bibitem{KilbasSrivastavaTrujillo2006}
A.~A. Kilbas, H.~M. Srivastava, J.~J. Trujillo, \emph{Theory and Applications
  of Fractional Differential Equations}. Vol. 204 of \emph{North-Holland
  Mathematics Studies}, Elsevier Science B.V., Amsterdam (2006).

\bibitem{KilbasTrujillo2001}
A.~A. Kilbas, J.~J. Trujillo, Differential equations of fractional order:
  methods, results and problems, {I}. \emph{Appl. Anal.} \textbf{78}, No 1-2
  (2001), 153--192.

\bibitem{Kiryakova1999}
V. Kiryakova, Multiindex {M}ittag-{L}effler functions, related
  {G}elfond-{L}eontiev operators and {L}aplace type integral transforms.
  \emph{Fract. Calc. Appl. Anal.} \textbf{2}, No~4 (1999), 445--462.

\bibitem{Kiryakova2000}
V.S. Kiryakova, Multiple (multiindex) {M}ittag-{L}effler functions and
  relations to generalized fractional calculus. \emph{J. Comput. Appl. Math.}
  \textbf{118}, No~1-2 (2000), 241--259.
  
 \bibitem{Kiryakova2010} 
 V. Kiryakova, The multi-index Mittag-Leffler functions as important class of special
 functions of fractional calculus. \emph{Computers and Math. with Appl.}
   {\bf 59}, No~5 (2010), 1885--1895.   

\bibitem{klafter2011first}
J.~Klafter, I.~M. Sokolov, \emph{{First Steps in Random Walks: From Tools to
  Applications}}. Oxford University Press (2011).

\bibitem{Kochubei2011}
A.~N. Kochubei, {General fractional calculus, evolution equations, and renewal
  processes}. \emph{Integr. Equat. Oper. Th.} \textbf{71} (2011), 583--600.

\bibitem{KonnoTamura2018}
H.~Konno, Y.~Tamura, {Stochastic modeling for neural spiking events based on
  fractional superstatistical Poisson process}. \emph{AIP Adv.} \textbf{8}
  (2018), 015118.

\bibitem{Krageloh2003}
A.~M. Kr\"{a}geloh, Two families of functions related to the fractional powers
  of generators of strongly continuous contraction semigroups. \emph{J. Math.
  Anal. Appl.} \textbf{283}, No~2 (2003), 459--467.

\bibitem{Laskin}
N.~Laskin, {Fractional Poisson process}. \emph{Comm. Nonlin. Sci. Numer.
  Simul.} \textbf{8} (2003), 201--213.

\bibitem{CC-exp1}
Z.~Lin, {On the FDTD formulations for biological tissues with Cole-Cole
  dispersion}. \emph{IEEE Microw. Wirel. Compon. Lett.} \textbf{20}, No~5
  (2010), 244--246.

\bibitem{Lubich1988a}
C.~Lubich, {Convolution quadrature and discretized operational calculus, {I}}.
  \emph{Numer. Math.} \textbf{52}, No~2 (1988), 129--145.

\bibitem{Lubich1988b}
C.~Lubich, {Convolution quadrature and discretized operational calculus, {II}}.
  \emph{Numer. Math.} \textbf{52}, No~4 (1988), 413--425.

\bibitem{Lubich2004}
C.~Lubich, {Convolution quadrature revisited}. \emph{BIT} \textbf{44}, No~3
  (2004), 503--514.
  
  \bibitem{Luchko1999} 
  Yu. Luchko, Operational method in fractional calculus. \emph{Fract. Calc. Appl. Anal.}  \textbf{2}, No 4 (1999), 463--489.

\bibitem{Mainardi2010}
F.~Mainardi, \emph{Fractional Calculus and Waves in Linear Viscoelasticity}.
  Imperial College Press, London (2010).

\bibitem{MainardiGarrappa2015_JCP}
F.~Mainardi, R.~Garrappa, On complete monotonicity of the {P}rabhakar function
  and non-{D}ebye relaxation in dielectrics. \emph{J. Comput. Phys.}
  \textbf{293} (2015), 70--80.

\bibitem{MainardiGorenflo2000}
F.~Mainardi, R.~Gorenflo, On {M}ittag-{L}effler-type functions in fractional
  evolution processes. \emph{J. Comput. Appl. Math.} \textbf{118}, No 1-2
  (2000), 283--299.

\bibitem{MainardiGorenflo2007}
F.~Mainardi, R.~Gorenflo, {Time-fractional derivatives in relaxation processes:
  a tutorial survey}. \emph{Fract. Calc. Appl. Anal.} \textbf{10}, No~3 (2007),
  269--308.

\bibitem{g3}
F.~Mainardi, R.~Gorenflo, E.~Scalas, {A renewal process of Mittag--Leffler
  type}. In: \emph{Thinking in Patterns}, 35--46, Word Scientific (2004).

\bibitem{MainardiSpada2011}
F.~Mainardi, G.~Spada, {Creep, relaxation and viscosity properties for basic
  fractional models in rheology}. \emph{Eur. Phys. J. Spec. Top.} \textbf{193},
  No~1 (2011), 133--160.

\bibitem{mathai2009h}
A.~M. Mathai, R.~K. Saxena, H.~J. Haubold, \emph{{The $H$-Function: Theory and
  Applications}}. Springer Science \& Business Media (2009).

\bibitem{meer}
M.~M. Meerschaert, E.~Nane, P.~Vellaisamy, {The fractional Poisson process and
  the inverse stable subordinator}. \emph{Electron. J. Probab.} \textbf{16},
  No~59 (2011), 1600--1620.

\bibitem{michelitsch2019continuous}
T.~M. Michelitsch, A.~P. Riascos, {Continuous time random walk and diffusion
  with generalized fractional Poisson process}. \emph{Physica A} (Online 31 Oct. 2019), 123294.

\bibitem{Michelitsch2019}
T.~M. Michelitsch, A.~P. Riascos, {Generalized fractional Poisson process and
  related stochastic dynamics}, arXiv: 1906.09704 (2019).
  
\bibitem{MillerRoss1993}
K.~S. Miller, B.~Ross, \emph{An Introduction to the Fractional Calculus and
  Fractional Differential Equations}. A Wiley-Intersci. Publ., John
  Wiley \& Sons, Inc., New York (1993).

\bibitem{Mittag-Leffler1902}
M.~G. Mittag-Leffler, Sur l'int\'{e}grale de {L}aplace-{A}bel. \emph{C. R.
  Acad. Sci. Paris (Ser. II)} \textbf{136} (1902), 937--939.

\bibitem{MittagLeffler1904}
M.~G. Mittag-Leffler, Sopra la funzione ${E}_{\alpha}(x)$. \emph{Rend. Accad.
  Lincei} \textbf{13}, No~5 (1904), 3--5.

\bibitem{montroll1965random}
E.~W. Montroll, G.~H. Weiss, {Random walks on lattices, II}. \emph{J. Math.
  Phys.} \textbf{6}, No~2 (1965), 167--181.

\bibitem{exp2}
K.~L. Ngai, A.~K. Jonscher, C.~T. White, {On the origin of the universal
  dielectric response in condensed matter}. \emph{Nature} \textbf{277}, No 5693
  (1979), 185--189.

\bibitem{NigmatullinOsokinSmith2003}
R.~Nigmatullin, S.~Osokin, G.~Smith, {The justified data-curve fitting
  approach: recognition of the new type of kinetic equations in fractional
  derivatives from analysis of raw dielectric data}. \emph{J. Phys. D}
  \textbf{36}, No~18 (2003), 2281--2294.

\bibitem{NigmatullinRyabov1997}
R.~Nigmatullin, Y.~Ryabov, {Cole--Davidson dielectric relaxation as a
  self-similar relaxation process}. \emph{Phys. Solid State} \textbf{39}, No~1
  (1997), 87--90.

\bibitem{NovikovWojciechowskiKomkovaThiel2005}
V.~Novikov, K.~Wojciechowski, O.~Komkova, T.~Thiel, {Anomalous relaxation in
  dielectrics. {E}quations with fractional derivatives}. \emph{Mater. Sci.
  Poland} \textbf{23}, No~4 (2005), 977--984.

\bibitem{paneva2013multi}
J.~Paneva-Konovska, {On the multi-index ($3m$-parametric) Mittag-Leffler
  functions, fractional calculus relations and series convergence}. \emph{Open
  Phys.} (\emph{Centr. Eur. J. Phys.}) \textbf{11}, No~10 (2013), 1164--1177.

\bibitem{paneva2014convergence}
J.~Paneva-Konovska, {Convergence of series in three parametric Mittag-Leffler
  functions}. \emph{Math. Slovaca} \textbf{64}, No~1 (2014), 73--84.

\bibitem{paneva2017overconvergence}
J.~Paneva-Konovska, {Overconvergence of series in generalized Mittag-Leffler
  functions}. \emph{Fract. Calc. Appl. Anal.} \textbf{20}, No~2 (2017),
  506--520.

\bibitem{Paris2010}
R.~B. Paris, Exponentially small expansions in the asymptotics of the {W}right
  function. \emph{J. Comput. Appl. Math.} \textbf{234}, No~2 (2010), 488--504.

\bibitem{Paris2019}
R.~B. Paris, Asymptotics of the special functions of fractional calculus. In:
  \emph{Handbook of Fractional Calculus with Applications}, {V}ol. 1, 297--325,
  De Gruyter, Berlin (2019).

\bibitem{Pipkin}
A.~C. Pipkin, \emph{Lectures on Viscoelasticity Theory}. Springer-Verlag
  (1972).

\bibitem{Podlubny1999}
I.~Podlubny, \emph{Fractional Differential Equations}. Vol. 198 of
  {Mathematics in Science and Engineering}. Academic Press Inc., San
  Diego, CA (1999).

\bibitem{Pogany2016}
T.~K. Pog\'any, Z.~Tomovski, {Probability distribution built by Prabhakar
  function. Related Tur\'an and Laguerre inequalities}. \emph{Integr. Transf. Spec. Funct.}
   \textbf{27}, No~10 (2016), 783--793.

\bibitem{scalas}
M.~Politi, T.~Kaizoji, E.~Scalas, {Full characterization of the fractional
  Poisson process}. \emph{EPL} \textbf{96} (2011), 20004.

\bibitem{Modern-stochastics-Polito}
F.~Polito, {Studies on generalized Yule models}. \emph{Mod. Stoch. Theory
  Appl.} \textbf{6} (2019), 41--55.

\bibitem{PolitoTomovski2016}
F.~Polito, {\v{Z}}.~Tomovski, Some properties of {P}rabhakar-type fractional
  calculus operators. \emph{Fract. Differ. Calc.} \textbf{6}, No~1 (2016),
  73--94.

\bibitem{Pollard1948}
H.~Pollard, The completely monotonic character of the {M}ittag-{L}effler
  function {$E_a(-x)$}. \emph{Bull. Amer. Math. Soc.} \textbf{54} (1948),
  1115--1116.

\bibitem{Prabhakar1971}
T.~R. Prabhakar, A singular integral equation with a generalized
  {M}ittag--{L}effler function in the kernel. \emph{Yokohama Math. J.}
  \textbf{19}, No~1 (1971), 7--15.

\bibitem{Repin}
O.~N. Repin, A.~I. Saichev, {Fractional Poisson law}. \emph{Radiophys. Quantum
  Electron.} \textbf{43}, No~9 (2000), 738--741.

\bibitem{RogosinMainardi-2014}
S.~Rogosin, F.~Mainardi, {George William Scott Blair -- The pioneer of
  factional calculus in rheology}. \emph{Commun. Appl. Ind. Math.} \textbf{6},
  No~1 (2014), e-481.

\bibitem{CC-exp2}
T.~Said, V.~V. Varadan, {Variation of Cole-Cole model parameters with the
  complex permittivity of biological tissues}. In: \emph{2009 IEEE MTT-S
  International Microwave Symposium Digest} (2009), 1445--1448.

\bibitem{Samko-book}
S.~G. Samko, A.~A. Kilbas, O.~I. Marichev, 
\emph{Fractional Integrals and Derivatives: Theory and Applications}. 
  Gordon and Breach Sci. Publ., Switzerland (1993).

\bibitem{Sandev-Diffusion7}
T.~Sandev, {Generalized Langevin equation and the Prabhakar derivative}.
  \emph{Mathematics} \textbf{5} (2017), \# 66.

\bibitem{Sandev-EPL}
T.~Sandev, I.~A, {Finite-velocity diffusion on a comb}. \emph{EPL}
  \textbf{124}, No~2 (2018), 20005.

\bibitem{Sandev-Diffusion1}
T.~Sandev, A.~Chechkin, H.~Kantz, R.~Metzler, {Diffusion and
  Fokker-Planck-Smoluchowski equations with generalized memory kernel}.
  \emph{Fract. Calc. Appl. Anal.} \textbf{18}, No~4 (2015), 1006--1038.

\bibitem{Sandev-Diffusion2}
T.~Sandev, A.~Chechkin, N.~Korabel, H.~Kantz, I.~Sokolov, R.~Metzler,
  {Distributed-order diffusion equations and multifractality: models and
  solutions}. \emph{Phys. Rev. E} \textbf{92} (2015), 042117.

\bibitem{Sandev-Diffusion5}
T.~Sandev, W.~Deng, P.~Xu, {Models for characterizing the transition among
  anomalous diffusions with different diffusion exponents}. \emph{J. Phys. A}
  \textbf{51}, No~40 (2018), 405002.

\bibitem{Sandev-Diffusion4}
T.~Sandev, R.~Metzler, A.~Chechkin, {From continuous time random walks to the
  generalized diffusion equation}. \emph{Fract. Calc. Appl. Anal.} \textbf{21},
  No~1 (2018), 10--28.

\bibitem{Sandev2019GeneralisedDA}
T.~Sandev, R.~Metzler, A.~V. Chechkin, Generalised diffusion and wave
  equations: Recent advances. In: \emph{Analytical Methods of Analysis and
  Differential Equations. AMADE-2018}. Cambridge Scientific Publishers (2019).

\bibitem{SandevTomovski2010}
T.~Sandev, {\v{Z}}.~Tomovski, Asymptotic behavior of a harmonic oscillator
  driven by a generalized {M}ittag-{L}effler noise. \emph{Phys. Scr.}
  \textbf{82}, No~6 (2010), \# 065001.

\bibitem{Sandev-Diffusion3}
T.~Sandev, {\v{Z}}.~Tomovski, {Langevin equation for a free particle driven by
  power law type of noises}. \emph{Phys. Lett. A} \textbf{378}, No~1--2 (2014),
  1--9.

\bibitem{SandevTomovski-Book}
T.~Sandev, {\v{Z}}.~Tomovski, \emph{Fractional Equations and Models: Theory
  and Applications}. Springer (2019).

\bibitem{Sandev-Diffusion6}
T.~Sandev, {\v{Z}}.~Tomovski, J.~L.~A. Dubbeldam, A.~Chechkin, {Generalized
  diffusion-wave equation with memory kernel}. \emph{J. Phys. A} \textbf{52},
  No~1 (2019), \# 015201.

\bibitem{SaxenaMathaiHaubold2006}
R.~K. Saxena, A.~M. Mathai, H.~J. Haubold, Reaction-diffusion systems and
  nonlinear waves. \emph{Astrophys. Space Sci.} \textbf{305}, No~3 (2006),
  297--303.

\bibitem{SaxenaGianni2011}
R.~K. Saxena, G.~Pagnini, Exact solutions of triple-order time-fractional
  differential equations for anomalous relaxation and diffusion, I: 
  The accelerating case. \emph{Physica A} \textbf{390}, No~4 (2011), 602--613.

\bibitem{SrivastavaTomovski2009}
H.~M. Srivastava, {\v{Z}}.~Tomovski, Fractional calculus with an integral
  operator containing a generalized {M}ittag-{L}effler function in the kernel.
  \emph{Appl. Math. Comput.} \textbf{211}, No~1 (2009), 198--210.

\bibitem{stanislavsky2018transient}
A.~Stanislavsky, A.~Weron, {Transient anomalous diffusion with Prabhakar-type
  memory}. \emph{J. Chem. Phys.} \textbf{149}, No~4 (2018), 044107.

\bibitem{StanislavskyWeron2012}
A.~Stanislavsky, K.~Weron, Numerical scheme for calculating of the fractional
  two-power relaxation laws in time-domain of measurements. \emph{Comput. Phys.
  Commun.} \textbf{183}, No~2 (2012), 320--323.

\bibitem{Stanislavsky2007}
A.~A. Stanislavsky, The stochastic nature of complexity evolution in the
  fractional systems. \emph{Chaos, Solitons $\&$ Fractals} \textbf{34}, No~1
  (2007), 51--61.

\bibitem{TomovskiHilferSrivastava2010}
{\v{Z}}.~Tomovski, R.~Hilfer, H.~M. Srivastava, Fractional and operational
  calculus with generalized fractional derivative operators and
  {M}ittag-{L}effler type functions. \emph{Integr. Transf. Spec. Funct.}
  \textbf{21}, No~11 (2010), 797--814.

\bibitem{TomovskiPoganySrivastava2014}
{\v{Z}}.~Tomovski, T.~Pog{\'a}ny, H.~M. Srivastava, Laplace type integral
  expressions for a certain three-parameter family of generalized
  {M}ittag-{L}effler functions with applications involving complete
  monotonicity. \emph{J. Franklin Inst.} \textbf{351}, No~12 (2014),
  5437--5454.

\bibitem{Uchaikin2003}
V.~Uchaikin, Relaxation processes and fractional differential equations.
  \emph{Internat. J. Theoret. Phys.} \textbf{42}, No~1 (2003), 121--134.

\bibitem{Vaughan1979}
W.~E. Vaughan, {Dielectric Relaxation}. \emph{Annu. Rev. Phys. Chem.}
  \textbf{30} (1979), 103--124.

\bibitem{WangCD2019}
C.~L. Wang, {Photocatalytic degradation as Davidson-Cole relaxation in time
  domain}. \emph{J. Adv. Dielectr.} \textbf{9}, No~1 (2019), \# 1950006.

\bibitem{WangHN2019}
C.~L. Wang, {Piezo-catalytic degradation of Havriliak-Negami type}. \emph{J.
  Adv. Dielectr.} \textbf{9}, No~3 (2019), \# 1950021.

\bibitem{WeidemanTrefethen2007}
J.~A.~C. Weideman, L.~N. Trefethen, Parabolic and hyperbolic contours for
  computing the {B}romwich integral. \emph{Math. Comp.} \textbf{76}, No 259
  (2007), 1341--1356.

\bibitem{WeronJurlewiczMagdziarz2005}
K.~Weron, A.~Jurlewicz, M.~Magdziarz, {H}avriliak-{N}egami response in the
  framework of the continuous-time random walk. \emph{Acta Phys. Pol. B}
  \textbf{36}, No~5 (2005), 1855--1868.

\bibitem{Wiman1905}
A.~Wiman, \"{U}ber den {F}undamentalsatz in der {T}eorie der {F}unktionen
  {$E^a(x)$}. \emph{Acta Math.} \textbf{29}, No~1 (1905), 191--201.

\bibitem{Wright1935}
E.~M. Wright, The asymptotic expansion of the generalised hypergeometric
  function. \emph{J. London Math. Soc.} \textbf{s1-10}, No~4 (1935), 286--293.

\bibitem{Wrigth1940_PTRSL}
E.~M. Wright, The asymptotic expansion of integral functions defined by
  {T}aylor series. \emph{Philos. Trans. Roy. Soc. London, Ser. A.} \textbf{238}
  (1940), 423--451.

\bibitem{Wright1940}
E.~M. Wright, The asymptotic expansion of the generalized hypergeometric
  function. \emph{Proc. London Math. Soc. (Ser. 2)} \textbf{46} (1940),
  389--408.

\end{thebibliography}
\end{document}